\documentclass[12pt]{article}

\usepackage{amssymb,amsmath}

\usepackage{pdfpages}

\newtheorem{theorem}{Theorem}[section]
\newtheorem{corollary}[theorem]{Corollary}

\newtheorem{proposition}[theorem]{Proposition}

\title{On Shilnikov's scenario with a homoclinic orbit in 3D - revised version}

\author{Hans-Otto Walther}

\begin{document}
	
\maketitle

\section{Introduction}

The purpose of the present account is to provide a detailed proof that complicated motion exists in {\it Shilnikov's scenario} \cite{S3}, which consists of a smooth vectorfield $V:\mathbb{R}^3\to\mathbb{R}^3$ with $V(0)=0$ so that the equation
\begin{equation}
x'(t)=V(x(t))
\end{equation}
has a homoclinic solution $h:\mathbb{R}\to\mathbb{R}^3\setminus\{0\}$ with $lim_{|t|\to\infty}h(t)=0$, and $DV(0)$ has eigenvalues $u>0$ and $\sigma\pm i\,\mu$, $\sigma<0<\mu$, with

\medskip

(H)$\quad\quad 0<\sigma+u.$

\medskip

Of course, smoothness of $V$ is to be specified in this statement. Shilnikov's result in \cite{S3} is that for $V$ analytic  there is a countable set of periodic orbits close to the homoclinic orbit $h(\mathbb{R})$.

\medskip

A related, slightly stronger statement about complicated motion is conjugacy with the shift $(s_j)_{-\infty}^{\infty}\mapsto(s_{j+1})_{-\infty}^{\infty}$
in two symbols $s_j\in\{0,1\}$, for a return map, 
which is given by intersections of solutions with a transversal to the homoclinic orbit. Work on verification of this property and of related ones can be found in the monographs  
\cite{GH,Wi1,Wi2,HSD} and in the references given there. A strong simplifying hypothesis is that $V$ is linear close to the equilibrium. For a detailed presentation under this hypothesis see \cite{HSD}. 

\medskip
	
Shilnikov-type results on complicated motion close to a homoclinic orbit approaching a stationary state  have also been obtained for semiflows in infinite-dimensional spaces, e. g. in \cite{W,LWW}
on delay differential equations which are linear close to the stationary state. These results are related not to \cite{S3} but to Shilnikov's work on vectorfields on $\mathbb{R}^4$ \cite{S4} which is analogous to \cite{S3} except for the presence of pairs of complex conjugate eigenvalues of $DV(0)$ in either halfplane.

\medskip

In the sequel we consider Shilnikov's scenario in $\mathbb{R}^3$ for $V$ twice continuously differentiable.
The hypothesis of second order differentiability, as opposed to minimal smoothness of merely continuous differentiability, serves  for the justification of the facilitating additional assumption that for $V$ (just once)  continuously differentiable both eigenspaces of $DV(0)$ are invariant under $V$ in a vicinity of $x=0$. The Appendix in Section 9 describes how the said invariance property can be achieved by means of a transformation which, however, reduces the order of differentiability by one.

\medskip

The result on complicated motion is stated in the final Theorem 8.2. Its proof uses the scaled form 
\begin{equation}
y'(t)=V_{\epsilon}(y(t)),\quad\epsilon>0,
\end{equation}
of Eq. (1), with $V_{\epsilon}:\mathbb{R}^3\to\mathbb{R}^3$ given by $V_{\epsilon}(x)=\frac{1}{\epsilon}V(\epsilon\,x)$. Eq. (2)  is equivalent to Eq.(1) as $y$ is a solution to Eq. (2) if and only if $x=\epsilon\,y$ is a solution to Eq. (1). 

\medskip

Following \cite{S3} we introduce a return map  in a transversal to the homoclinic orbit,  with domain on one side of the local stable manifold $W^s_{loc}$ at $y=0$,   ''above" or ''below" $W^s_{loc}$, compare Figure 1 on page 3. 

\begin{figure}
\includegraphics[page=1,scale=0.7]{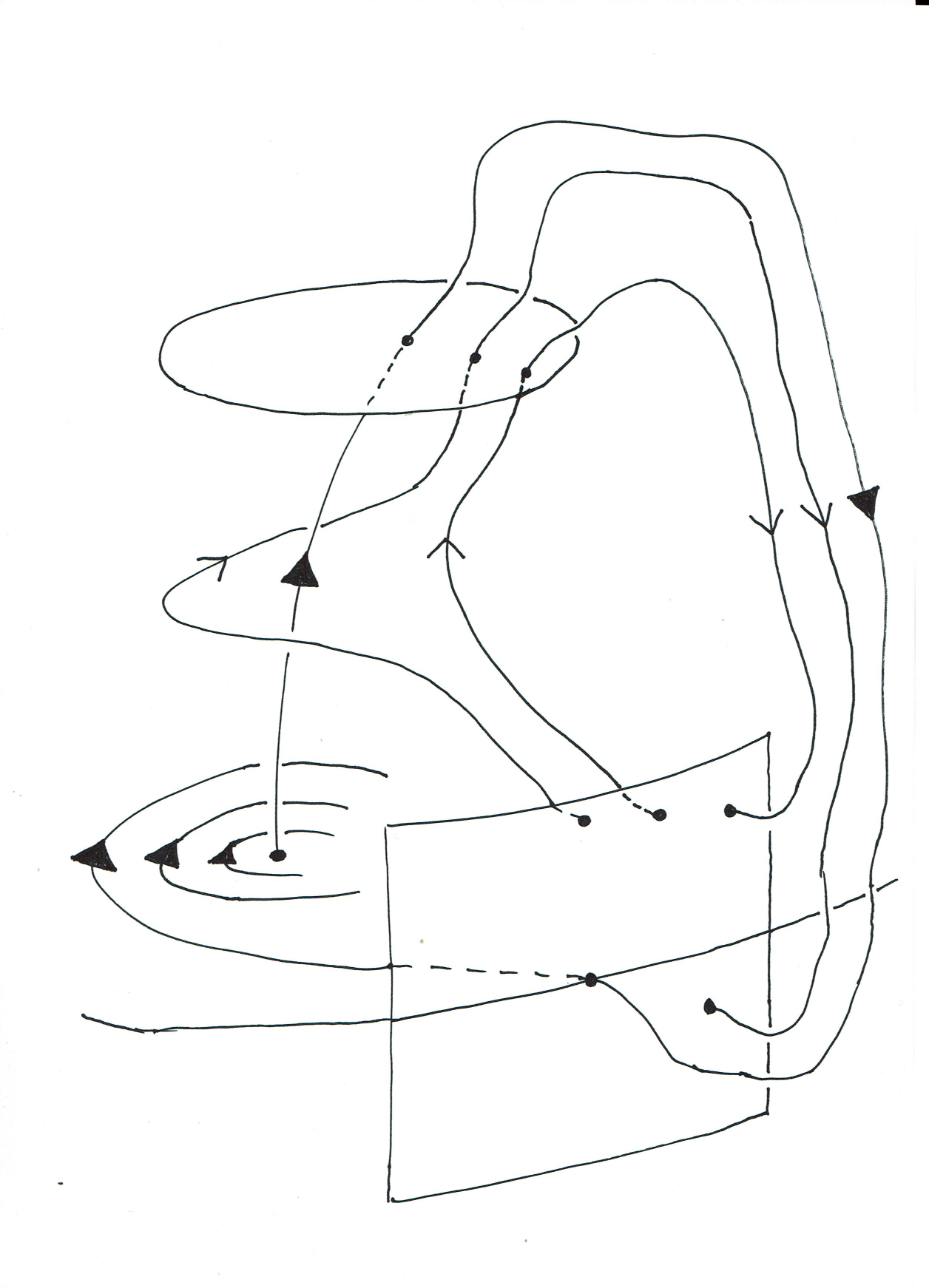}
\caption{Actions of the return map}
\end{figure}

The return map is transformed to a map in the plane. For $\epsilon>0$ sufficiently small we find disjoint sets $M_0,M_1$ so that for each symbol sequence $(s_j)_{j=-\infty}^{\infty}$
there exist trajectories $(x_j)_{j=-\infty}^{\infty}$ of the transformed return map with
$x_j\in M_{s_j}$ for all integers $j$.

\medskip

The result is achieved by an examination of the expansion of curves under the return map, and by the inspection of  intersections of the expanded curves with the domain of the return map.

\medskip

The reader may observe that the proof of Theorem 8.2 actually yields a countable family of sets of complicated trajectories of the return map.

\medskip

No attempt is made in the present paper to verify covering relations like for example a Smale horseshoe, and neither existence of periodic orbits close to the homoclinic nor conjugacy of the return map on some invariant set with the shift in two symbols are touched upon.

\medskip

{\bf Notation, preliminaries}. A forward trajectory of a map $f:M\supset\,dom\to M$ is a sequence $(x_j)_0^{\infty}$ in $dom$ with with $x_{j+1}=f(x_j)$ for all integers $j\ge0$. 
Entire trajectories are defined analogously, with all integers as indices.
For a vectorspace $X$, $x\in X$ and $M\subset X$ we set $x\pm M=\{y\in X:y\mp x\in M\}$. Similarly, for $A\subset\mathbb{R}$ and $x\in X$, $Ax=\{y\in X:\mbox{For some}\,\, a\in A, y=ax\}$. The interior, the boundary, and the closure of a subset of a topological space are denoted by $int\,M, \partial \,M,$ and $cl\,M,$ respectively. Components of vectors in Euclidean spaces $\mathbb{R}^n$ are indicated by lower indices.The inner product on $\mathbb{R}^n$ is written as $<x,y>=\sum_{i=1}^nx_iy_i$, and we use the Euclidean norms given by $|x|=\sqrt{<x,x>}$. The vectors of the canonical orthonormal basis on $\mathbb{R}^n$ are denoted by $e_j$, $j=1,\ldots,n$, $e_{jj}=1$ and $e_{jk}=0$ for $j\neq k$. In $\mathbb{R}^3$ we write $L=\mathbb{R}e_1\oplus\mathbb{R}e_2$ and $U=\mathbb{R}e_3$. The  associated projections onto $L$ and onto $U$ are  $P_L:\mathbb{R}^3\to\mathbb{R}^3$ and $P_U:\mathbb{R}^3\to\mathbb{R}^3$, respectively. For every $x\in\mathbb{R}^3$, $|P_Ux|^2+|P_Lx|^2=|x|^2$, and each of the projections has norm $1$  in $L_c(\mathbb{R}^3,\mathbb{R}^3)$.  
For a function $f:\mathbb{R}^n\supset dom\to\mathbb{R}$ and for $x\in dom$ the derivatives and partial derivatives are linear maps $Df(x)$ and $D_jf(x)$, respectively. Partial derivatives as numbers satisfy $\partial_jf(x)=D_jf(x)1=Df(x)e_j$. The tangent space $T_xM$ of a continuously differentiable submanifold $M$ of $\mathbb{R}^n$ at $x\in M$ is the set of tangent vectors $v=c'(0)=Dc(0)1$ of continuously differentiable curves $c:I\to\mathbb{R}^n$ with $I$ an interval, $c(I)\subset M$, $0\in I$, $c(0)=x$. For a continuously differentiable map $f:M\to N$, $N$ a continuously differentiable submanifold of  $\mathbb{R}^n$, the derivative at $x\in M$ is the linear continuous map $T_xf:T_xM\to T_{f(x)}N$ given by $T_xf(v)=(f\circ c)'(0)$, for $v$ and $c$ as before. The flow
${\mathcal F}$ generated by a vectorfield ${\mathcal V}:\mathbb{R}^n\supset{\mathcal U}\to\mathcal{R}^n$ which is locally Lipschitz continuous is the map $\mathbb{R}\times\mathbb{R}^n\supset dom_{{\mathcal F}}\to\mathbb{R}^n$ 
which is given by $(t,x)\in\,dom_{{\mathcal F}}$ if and only if $t$ belongs to the domain of the maximal solution $y:I_x\to\mathbb{R}^n$ of the differential equation $x'(t)={\mathcal V}(x(t))$ with initial value $y(0)=x$, and ${\mathcal F}(t,x)=y(t)$ in this case. ${\mathcal F}$ is of the same order of differentiability as ${\mathcal V}$. A subset $M\subset U$ is invariant under ${\mathcal F}$ if $x\in M$ implies ${\mathcal F}(t,x)\in M$ for all $t\in\mathbb{R}$ with $(t,x)\in \,dom_{{\mathcal F}}$. For a further subset $N\subset U$ the set $M$ is called invariant under ${{\mathcal F}}$ in $N$ if for every $x\in M\cap N$ and for every interval $I\ni 0$ with ${{\mathcal F}}(I\times\{x\})\subset N$ we have ${{\mathcal F}}(I\times\{x\})\subset M$.

\section{A simpler setting, and scaled vectorfields}

We assume from here on that we are in Shilnikov's scenario with $V$ continuously differentiable, and that in addition $V$ has the subsequent properties (A1) and (A2). In the Appendix (Section 9) it is explained how (A1) and (A2) can be achieved if we start from Shilnikov's scenario for a vectorfield which is twice continuously differentiable. 

\medskip

(A1) For all $x\in\mathbb{R}^3$, $DV(0)x=Ax$ with a matrix
\begin{equation*}
A=
\left(
\begin{array}{ccc}
	\sigma & \mu & 0\\
	-\mu & \sigma & 0\\
	0 & 0 & u
\end{array}
\right).
\end{equation*}
Then $AL=L$ and $AU=U$ are the (realified) eigenspaces associated with the eigenvalues $\sigma\pm i\,\mu$ and $u$, respectively. 

\medskip

(A2) There  exists $r_V>0$ with 
\begin{eqnarray}
P_UV(x)=0\quad\mbox{for}\quad x\in L,\quad|x|\le r_V,\\
P_LV(x)=0\quad\mbox{for}\quad x\in U,\quad|x|\le r_V,
\end{eqnarray}  
and there are
reals $t_U<t_L$ with
\begin{eqnarray}
h(t)\in U\quad\mbox{for}\quad t\le t_U,\\
h(t)\in L\quad\mbox{for}\quad t\ge t_L.
\end{eqnarray}

In order to exploit the hypothesis (H) on eigenvalues we need to study solutions of Eq.(1) in small neighbourhoods of the origin where $V$ is appropriately close to its linearization $DV(0)$.  
Alternatively one can study solutions of the equivalent Eqs. (2) for small $\epsilon>0$ in a fixed neighbourhood of the origin.

We prefer the latter and follow solutions of the Eqs. (2) inside the compact truncated solid cylinder
$$
B_1=\{x\in\mathbb{R}^n:|P_Ux|\le1,|P_Lx|\le1\},
$$
Let $F:\mathbb{R}^3\supset\,dom\to\mathbb{R}^3$ denote the flow generated by Eq.(1). 

\medskip

The scaled vectorfields
$V_{\epsilon}:\mathbb{R}^3\to\mathbb{R}^3$, $\epsilon>0$, are continuously differentiable and satisfy $V_{\epsilon}(0)=0$. For every $\epsilon>0$, 

\medskip

- the solution $h_{\epsilon}=\frac{1}{\epsilon}h$ of Eq. (2) satisfies $\lim_{|t|\to\infty}h_{\epsilon}(t)=0$,

\medskip

- we have
$$
DV_{\epsilon}(0)=DV(0)
$$
since $DV_{\epsilon}(z)y=\frac{1}{\epsilon}DV(\epsilon z)\epsilon y=DV(\epsilon z)y$ for $\epsilon>0,z\in\mathbb{R}^3,y\in\mathbb{R}^3$,

\medskip

- the flow $F_{\epsilon}:\mathbb{R}\times\mathbb{R}^3\supset\Omega_{\epsilon}\to\mathbb{R}^3$ generated by Eq. (2) is given by
$(t,x_0)\in\Omega_{\epsilon}$ if and only if $(t,\epsilon\,x_0)\in\Omega$, and $F_{\epsilon}(t,x_0)=\frac{1}{\epsilon}F(t,\epsilon\,x_0)$.

\medskip

The following proposition states precisely that in $B_1$ the linear spaces $L$ and $U$ are invariant under scaled vectorfields $V_{\epsilon}$.

\begin{proposition}
For $0<\epsilon\le r_V$ and for all $x\in B_1$,
$$
P_UV_{\epsilon}(P_Lx)=0\quad\mbox{and}\quad P_LV_{\epsilon}(P_Ux)=0.
$$
\end{proposition}
{\bf Proof.}
For $0<\epsilon\le r_V$ and $x\in B_1$, we have $P_Lx\in L\cap B_1$ and $|\epsilon P_Lx|=\epsilon |P_Lx|\le\epsilon\le r_V$. Using Eq. (3) we get $P_UV(\epsilon P_Lx)=0$, which in turn yields $0=P_U\frac{1}{\epsilon}V(\epsilon P_Lx)=P_UV_{\epsilon}(P_Lx)$. The other equation of the assertion is shown analogously.
$\Box$

\medskip

Proposition 2.1 can be rephrased as  
 $$
 V_{\epsilon}(L\cap B_1)\subset L\quad\mbox{and}\quad V_{\epsilon}(U\cap B_1)\subset U\quad\mbox{for}\quad0<\epsilon\le r_v.
 $$
The next proposition expresses closeness of scaled vectorfields to their linearization at the origin in terms of components in $U$ and $L$.

\begin{proposition}
For every $\eta>0$ there exists $\epsilon(\eta)\in(0,r_V]$ such that for all $\epsilon\in(0,\epsilon(\eta))$,
\begin{eqnarray}
|DV(\epsilon y)-A| & \le & \eta\,\,\mbox{for all}\,\,y\in B_1,\\
|P_U(V_{\epsilon}(x)-Ax)| & \le & \eta|P_Ux|\,\,\mbox{for all}\,\,x\in B_1,\\
|P_L(V_{\epsilon}(x)-Ax)| & \le & \eta|P_Lx|\,\,\mbox{for all}\,\,x\in B_1.
\end{eqnarray}
\end{proposition}

{\bf Proof.} Let $\eta>0$ be given. By continuity of $DV$ at $0$ there exists $\epsilon(\eta)\in (0,r_V]$ so that for all $\epsilon\in(0,\epsilon(\eta))$ and all $y \in B_1$,
$$
|DV(\epsilon y)-A|\le\eta,
$$
which is Eq. (7). Proof of Eq. (8):
\begin{eqnarray*}
|P_U(V_{\epsilon}(x)-Ax)| & = & |P_U(V_{\epsilon}(x)-Ax)-(P_U(V_{\epsilon}(P_Lx)-AP_Lx))|\\
& & \mbox{(the subtracted term is zero, by Proposition 2.1 and}\quad P_UAP_L=0)\\
& = & \left|\int_0^1P_U(DV_{\epsilon}(y_t))-A)[x-P_Lx]dt\right|\\
& & \mbox{where}\quad y_t=P_Lx+t(x-P_Lx)=P_Lx+tP_Ux\in B_1\\
& = & \left|\int_0^1P_U\left\{\frac{1}{\epsilon}DV(\epsilon y_t)\circ(\epsilon\cdot id)-A\right\}[P_Ux]dt\right|\\
& = & \left|\int_0^1P_U\{DV(\epsilon y_t)-A\}[P_Ux]dt\right|\le\max_{0\le t\le1}|P_U\{DV(\epsilon y_t)-A\}||P_Ux|\\
& \le & \max_{y\in B_1}|DV(\epsilon y)-A||P_Ux|\le\eta|P_Ux|\quad\mbox{(with}\quad|P_U|=1).
\end{eqnarray*}Eq. (9) is shown analogously. $\Box$

\medskip

Incidentally, notice that conversely the inequalities (8) and (9) for $x\in B_1$ imply that $V_{\epsilon}(U\cap B_1)\subset U$ and $V_{\epsilon}(L\cap B_1)\subset L$.

\begin{proposition}
There exists $\epsilon_M>0$ so that for every $\epsilon\in(0,\epsilon_M)$ there are $t_{E,\epsilon}<t_{I,\epsilon}$ with 
$$
h_{\epsilon}(t_{E,\epsilon})\in\{-e_3,e_3\}\quad\mbox{and}\quad h_{\epsilon}(t_{I,\epsilon})\in L,\quad|h_{\epsilon}(t_{I,\epsilon})|=1.
$$
Either $h_{\epsilon}(t)\in(-\infty,0)e_3$ for all $t\le t_{E,\epsilon}$, or	
\begin{equation}
h_{\epsilon}(t)\in(0,\infty)e_3\quad\mbox{for all}\quad t\le t_{E,\epsilon}.
\end{equation}
For all $t\ge t_{I,\epsilon}$, $h_{\epsilon}(t)\in L$.
\end{proposition}

{\bf Proof.} 1. For all $t\le t_U$, $h(t)\in U\setminus\{0\}=(-\infty,0)e_3\cup(0,\infty)e_3$. Hence either
$h(t)\in(-\infty,0)e_3$ for all $t\le t_U$, or $h(t)\in(0,\infty)e_3$ for all $t\le t_U$.

\medskip

2. On $t_{E,\epsilon}$ for $\epsilon\in(0,|h(t_U)|)$. In case $h(t)\in(0,\infty)e_3$ for all $t\le t_U$ the relation  $\lim_{t\to-\infty}h(t)=0$ shows that for each $\epsilon\in(0,|h(t_U)|)$ there exists $t_{E,\epsilon}\le t_U$ with $|h(t_{E,\epsilon})|=\epsilon$. It follows that $h(t_{E,\epsilon})=\epsilon e_3$. For the same $\epsilon<|h(t_U)|$ and for all $t\le t_{E,\epsilon}\le t_U$
we get $h_{\epsilon}(t)=\frac{1}{\epsilon}h(t)\in(0,\infty)e_3$ and in particular, $h_{\epsilon}(t_{E,\epsilon})=\frac{1}{\epsilon}h(t_{E,\epsilon})=e_3$.
Analogously for the second case.

\medskip

3. On $t_{I,\epsilon}$ for $\epsilon\in(0,|h(t_L)|)$. The relations $h(t)\in L$ for $t\ge t_L$ and $\lim_{t\to\infty}h(t)=0$ show that for each $\epsilon\in(0,|h(t_L)|)$ there exists $t_{I,\epsilon}\ge t_L\ge t_U\ge t_{E,\epsilon}$  with $|h(t_{I,\epsilon})|=\epsilon$. Hence $|h_{\epsilon}(t_{I,\epsilon})|=1$. For $t\ge t_{I,\epsilon}\ge t_L$,
$h_{\epsilon}(t)=\frac{1}{\epsilon}h(t)\in\frac{1}{\epsilon}L= L$. 

\medskip

4. Set $\epsilon_M=\min\{|h(t_U)|,|h(t_L)|\}$. $\Box$

\medskip

In the sequel we focus on case (10) of Proposition 2.3 - the other case is analogous.

\section{Transversality, projected solution curves, polar coordinates}

 We want to describe the behavior of solutions to Eq. (2)  close to the homoclinic loop $h_{\epsilon}(\mathbb{R})\cup\{0\}$  for small $\epsilon\in(0,\epsilon_M)$. This will be done in terms of a return map which follows the solutions from a subset of the cylinder
$$
M_I=\{x\in\mathbb{R}^3:|P_Lx|=1\}
$$
above the plane $L$ until their return to $M_I$. The return map will be obtained as a composition of an {\it inner map}, which follows solutions until they reach the plane
$$
M_E=e_3+L,
$$
parallel to $L$ (compare Figure 1 on page 3), with an {\it exterior map} following solutions from a subset of $M_E$ until they reach $M_I$. The construction and analysis of the inner and exterior maps require preparations.

\medskip

We begin with transversality of the vectorfield where the homoclinic solution intersects with $M_E$ and with $M_I$, which is at $t=t_{E,\epsilon}$ and at $t=t_{I,\epsilon}$, respectively. The tangent spaces of $M_E$ are all equal to $L$ while for every $x\in M_I$  we have
$$
T_xM_I=x^{\perp}\oplus U.
$$
where $x^{\perp}\in L$ is orthogonal to $P_Lx$ in $L$,
$x^{\perp}=\left(\begin{array}{c}-x_2\\x_1\\0\end{array}\right)$. 

\begin{proposition}
(i) There exists $\epsilon_E>0$ so that for all $\epsilon\in(0,\epsilon_E)$, $V_{\epsilon}(e_3)\notin T_{e_3}M_E$.

\medskip

(ii)  There exists $\epsilon_I>0$ so that for all $\epsilon\in(0,\epsilon_I)$ and for all $x\in M_I$ with $|x_3|<1$,
$V_{\epsilon}(x)\notin T_xM_I$.
\end{proposition}

{\bf Proof.} 1. On (i). From Proposition 2.2 for $0<\epsilon<\epsilon(u/2)$ we get
$$
|P_UV_{\epsilon}(e_3)|\ge|P_UAe_3|-|P_U(V_{\epsilon}(e_3)-Ae_3)|
=u-|P_U(V_{\epsilon}(e_3)-Ae_3)|
\ge u-\frac{u}{2}|P_Ue_3|=\frac{u}{2}>0.
$$
Hence $P_UV_{\epsilon}(e_3)\notin L=T_{e_3}M_E$.

\medskip

2. On (ii). For $x\in M_I$ with $|x_3|<1$, $x\in B_1$.  $P_Lx\neq0$ and $\left(\begin{array}{c}-x_2 \\ x_1 \\0\end{array}\right)\neq0$ are orthogonal. Using Proposition 2.2 for $0<\epsilon<\epsilon(\frac{|\sigma|}{2})$  we get
\begin{eqnarray*}
|\big<P_LV_{\epsilon}(x),P_Lx\big>| & \ge & |\big<P_LAx),P_Lx\big>| - |\big<P_L(V_{\epsilon}(x)-Ax),P_Lx\big>|\\
& \ge & \bigg|\bigg<\left(\begin{array}{c} \sigma x_1+\mu x_2 \\-\mu x_1+\sigma x_2 \\0\end{array}\right),\left(\begin{array}{c} x_1 \\ x_2 \\0\end{array}\right)\bigg>\bigg|-\frac{|\sigma|}{2}|P_Lx|^2\\
& = & (-\sigma)|P_Lx|^2-\frac{|\sigma|}{2}|P_Lx|^2>0,
\end{eqnarray*}
which yields $P_LV_{\epsilon}(x)\notin\mathbb{R}\left(\begin{array}{c}-x_2 \\ x_1 \\0\end{array}\right)$. Now the assertion follows easily. $\Box$

\begin{corollary}
For $0<\epsilon<\min\{\epsilon_M,\epsilon(u/2),\epsilon(|\sigma|/2)\}$,
$$
h_{\epsilon}'(t_{E,\epsilon})=V_{\epsilon}(h_{\epsilon}(t_{E,\epsilon}))\notin T_{h_{\epsilon}(t_{E,\epsilon})}M_E\quad\mbox{and}\quad
h_{\epsilon}'(t_{I,\epsilon})=V_{\epsilon}(h_{\epsilon}(t_{I,\epsilon}))\notin T_{h_{\epsilon}(t_{I,\epsilon})}M_I.
$$
\end{corollary}

\medskip

In order to formulate differential equations for projections of solutions of Eq. (2) into $L$ and $U$ we introduce
$R_{\epsilon,L}:\mathbb{R}^3\to\mathbb{R}^3$ and $R_{\epsilon,U}:\mathbb{R}^3\to\mathbb{R}^3$ by 
$$
R_{\epsilon,L}(x)=P_L(V_{\epsilon}(x)-Ax)\in L\quad\mbox{and}\quad R_{\epsilon,U}(x)=P_U(V_{\epsilon}(x)-Ax)\in U.
$$
\begin{corollary}
Assume 
$$
\mbox{H}_1(\eta,\epsilon)\qquad 0<\eta\quad\mbox{and}\quad0<\epsilon<\min\{\epsilon(\eta),\epsilon_M,\epsilon_I,\epsilon_E,\epsilon(u/2),\epsilon(-\sigma/2)\}.
$$
For every solution $y:I\to\mathbb{R}^3$ of Eq. (2) and $t\in I$ with $y(t)\in B_1$,
\begin{eqnarray*}
(P_Ly)'(t) & = & AP_Ly(t)+R_{\epsilon,L}(y(t))\quad\mbox{with}\quad|R_{\epsilon,L}(y(t))|\le\eta|P_Ly(t)|,\\
(P_Uy)'(t) & = & AP_Uy(t)+R_{\epsilon,U}(y(t))\quad\mbox{with}\quad|R_{\epsilon,U}(y(t))|\le\eta|P_Uy(t)|.
\end{eqnarray*}
\end{corollary}

\medskip

We turn to the position of projections $P_Ly$ of solutions $y$ to Eq. (2) in terms of polar coordinates. For the statement of the next result it is convenient to introduce the maps ${\mathcal A}_{\epsilon}:\{x\in\mathbb{R}^3: x_1^2+x_2^2>0\}\to\mathbb{R}^3$ and ${\mathcal B}_{\epsilon}:\{x\in\mathbb{R}^3: x_1^2+x_2^2>0\}\to\mathbb{R}^3$ given by
$$
{\mathcal A}_{\epsilon}(x)=\bigg<\left(
\begin{array}{c}
R_{\epsilon,L,1}(x)\\
R_{\epsilon,L,2}(x)
\end{array}
\right)
,
\frac{1}{x_1^2+x_2^2}	
\left(
\begin{array}{c}
x_1\\
x_2
\end{array}
\right)
\bigg>\quad\mbox{and}\quad {\mathcal B}_{\epsilon}(x)=\bigg<\left(
\begin{array}{c}
R_{\epsilon,L,1}(x)\\
R_{\epsilon,L,2}(x)
\end{array}
\right)
,
\frac{1}{x_1^2+x_2^2}	
\left(
\begin{array}{c}
-x_2\\
x_1
\end{array}
\right)
\bigg>.
$$

\begin{proposition}
 Assume H$_1(\eta,\epsilon)$. Let $I\subset\mathbb{R}$ be an interval and consider a solution $y:I\to\mathbb{R}^3$ of Eq. (2).

(i)  For every $t\in I$,
$$
y_3'(t)=u\,y_3(t)+R_{\epsilon,U,3}(y(t))
$$
and
$$
\left(
\begin{array}{c}
y_1'(t)\\
y_2'(t)
\end{array}
\right)	
 = 
\left( 
\begin{array}{cc}
	\sigma & \mu\\
	-\mu & \sigma
\end{array}
\right)	
\left(
\begin{array}{c}
	y_1(t)\\
	y_2(t)
\end{array}
\right)
+
\left(
\begin{array}{c}
R_{\epsilon,L,1}(y(t))\\
R_{\epsilon,L,2}(y(t))
\end{array}
\right)
.
$$
In case $y(t)\in B_1$,
$$
|(R_{\epsilon,U,3}(y(t))|\le\eta\,y_3(t)\quad\mbox{and}\quad
\left|\left(
\begin{array}{c}
	R_{\epsilon,L,1}(y(t))\\
	R_{\epsilon,L,2}(y(t))
\end{array}
\right)
\right|
\le \eta
\left|\left(
\begin{array}{c}
	y_1(t)\\
	y_2(t)
\end{array}
\right)
\right|.
$$

(ii) Assume $0\in I$, and $(y_1(t))^2+(y_2(t))^2>0$
for all $t\in I$. Let $\psi\in\mathbb{R}$ be given such that 
$$
\left(\begin{array}{c}y_1(0)\\y_2(0)\end{array}\right)=\sqrt{(y_1(0))^2+(y_2(0))^2}\left(\begin{array}{c}\cos(\psi)\\\sin(\psi)\end{array}\right).
$$ 
Then $y$ and the continuously differentiable functions $r:I\to(0,\infty)$ and $\phi:I\to\mathbb{R}$ given by
$$
r(t)=\sqrt{(y_1(t))^2+(y_2(t))^2}\quad\mbox{and}\quad\phi(t)=
\psi-\sigma\,t+\int_0^t{\mathcal B}_{\epsilon}(y(s))ds
$$
satisfy 
\begin{eqnarray}
y'(t) & = & V_{\epsilon}(y(t)),\nonumber\\	
r'(t) & = & \sigma\,r(t)+{\mathcal A}_{\epsilon}(y(t))r(t),\\	
\phi'(t) & = & -\mu+{\mathcal B}_{\epsilon}(y(t))
\end{eqnarray}
for every $t\in I$, with
$$
\left(
\begin{array}{c}
y_1(t)\\
y_2(t)
\end{array}
\right)
=r(t)	
\left(
\begin{array}{c}
\cos(\phi(t))\\
\sin(\phi(t))
\end{array}
\right)
\quad\mbox{for all}\quad t\in I.
$$

(iii) If $y(I)\subset B_1$ and $(y_1(t_0))^2+(y_2(t_0))^2>0$
for some $t_0\in I$ then $(y_1(t))^2+(y_2(t))^2>0$
for all $t\in I$.
\end{proposition}

Compare Figure 2 on page 9.

\medskip

\begin{figure}
\includegraphics[page=1,scale=0.7]{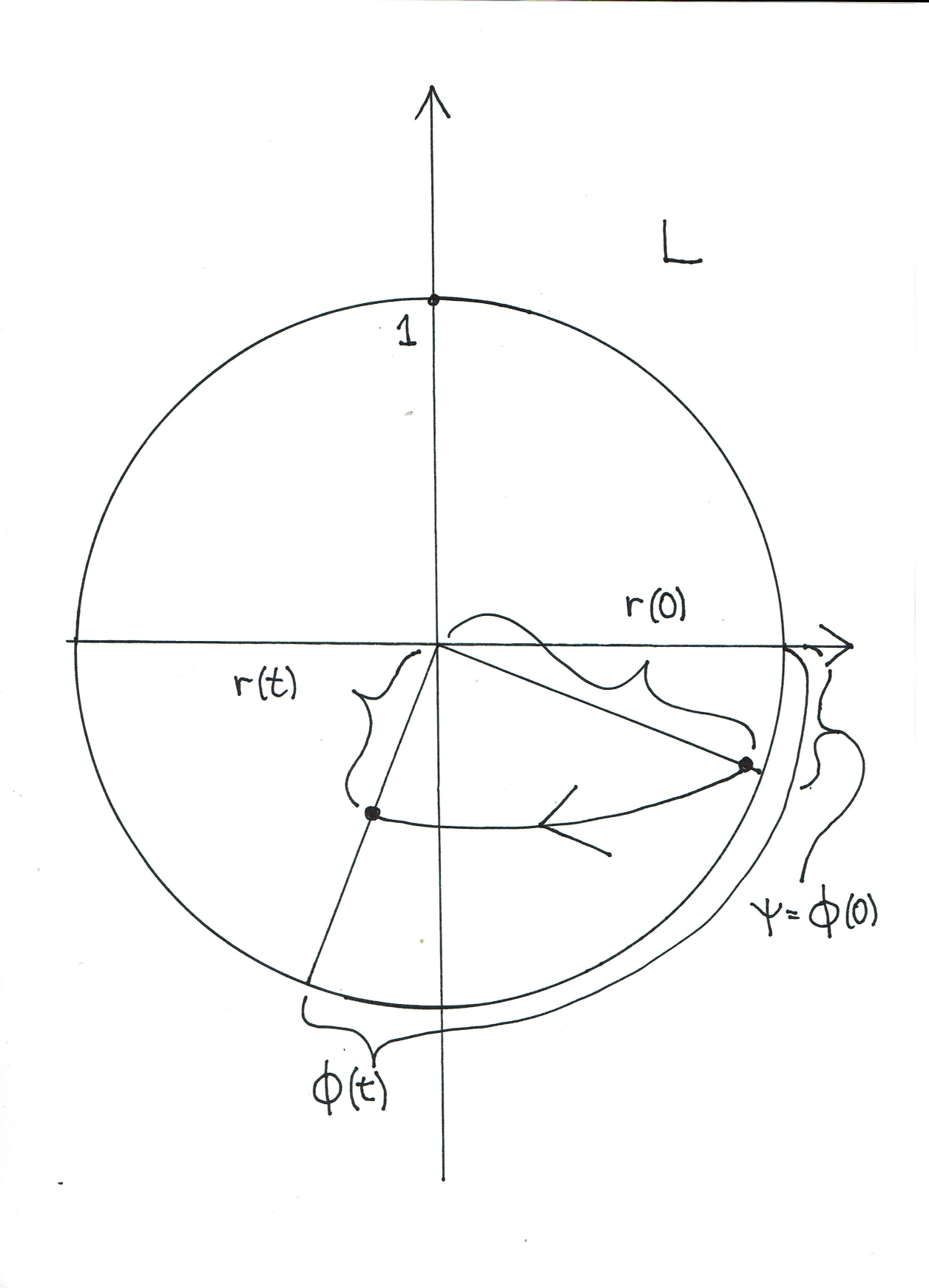}
\caption{The polar coordinates along the projected solution}
\end{figure}

{\bf Proof.} 1. Assertion (i) is a consequence of Eq. (7) in combination with Corollary 3.3.

\medskip

2. Proof of assertion (ii). 

\medskip

2.1. Proof of Eq. (11). Differentiation of $r^2$ in combination with the differential equation from assertion (i) yields
\begin{eqnarray*}
2r(t)r'(t) & = & 2\bigg<\left(\begin{array}{c}
	y_1(t)\\
	y_2(t)
\end{array}
\right),\left(\begin{array}{c}
	y_1'(t)\\
	y_2'(t)
\end{array}
\right)\bigg>\\
& = & 2\,\bigg<
\left(
\begin{array}{c}
	y_1(t)\\
	y_2(t)
\end{array}
\right),\left( 
\begin{array}{cc}
	\sigma & \mu\\
	-\mu & \sigma
\end{array}
\right)	
\left(
\begin{array}{c}
	y_1(t)\\
	y_2(t)
\end{array}
\right)\bigg>
+
2\,\bigg<
\left(
\begin{array}{c}
	R_{\epsilon,L,1}(y(t))\\
	R_{\epsilon,L,2}(y(t))
\end{array}
\right),
\left(
\begin{array}{c}
y_1(t)\\
y_2(t))\end{array}
\right)
\bigg>	\\
& = & 2(
y_1(t)[\sigma\,y_1(t)+\mu\,y_2(t)]+y_2(t)[-\mu\,y_1(t)+\sigma\,y_2(t)])+2\,{\mathcal A}_{\epsilon}(y(t))(r(t))^2)\\
& = & 2\,\sigma(r(t))^2+2\,{\mathcal A}_{\epsilon}(t)(r(t))^2.
\end{eqnarray*}
Divide by $2\,r(t)>0$.

\medskip

2.2. Differentiation of $\phi$ yields Eq. (12).

\medskip

2.3. Decomposition of the right hand side of the differential equation for $\left(\begin{array}{c}y_1\\y_2\end{array}\right)$  from assertion (i) with respect to the orthonormal bases 
$$
\left\{\frac{1}{r(t)}\left(\begin{array}{c}y_1(t)\\y_2(t)\end{array}\right),\frac{1}{r(t)}\left(\begin{array}{c}-y_2(t)\\y_1(t)\end{array}\right)\right\},\quad t\in I,
$$
yields that $\left(
\begin{array}{c}
y_1\\
y_2
\end{array}
\right)$ is a solution of the nonautonomous linear system
$$
\left(
\begin{array}{c}
y_1'(t)\\
y_2'(t)
\end{array}
\right)	
= 
\left( 
\begin{array}{cc}
\sigma & \mu\\
-\mu & \sigma
\end{array}
\right)	
\left(
\begin{array}{c}
y_1(t)\\
y_2(t)
\end{array}
\right)
+{\mathcal A}_{\epsilon}(y(t))
\left(
\begin{array}{c}
y_1(t)\\
y_2(t)
\end{array} 
\right)
(t)+{\mathcal B}_{\epsilon}(y(t))\left(
\begin{array}{c}
-y_2(t)\\
y_1(t)
\end{array}
\right).
$$

2.4. For the continuously differentiable map
$$
w:I\ni t\mapsto r(t)\left(\begin{array}{c}\cos(\phi(t))\\\sin(\phi(t))\end{array}\right)\in\mathbb{R}^2.
$$
we compute
\begin{eqnarray*}
w'(t) & = & r'(t)\left(\begin{array}{c}\cos(\phi(t))\\ \sin(\phi(t))\end{array}\right)+r(t)\phi'(t)\left(\begin{array}{c}-\sin(\phi(t))\\ \cos(\phi(t))\end{array}\right)\\
& = & [\sigma\,r(t)+{\mathcal A}_{\epsilon}(y(t))r(t)]\left(\begin{array}{c}\cos(\phi(t))\\ \sin(\phi(t))\end{array}\right)+r(t)[-\sigma+{\mathcal B}_{\epsilon}(y(t))]\left(\begin{array}{c}-\sin(\phi(t))\\ \cos(\phi(t))\end{array}\right)\\
& & \mbox{(with Eqs. (11,12))}\\
& = & \sigma\,w(t)-\mu\left(\begin{array}{c}-w_2(t)\\w_1(t)\end{array}\right)+{\mathcal A}_{\epsilon}(y(t))w(t)+{\mathcal B}_{\epsilon}(y(t))\left(\begin{array}{c}-w_2(t)\\w_1(t)\end{array}\right)\\
& = & \left(\begin{array}{c}\sigma\,w_1(t)+\mu\,w_2(t)\\-\mu\,w_1(t)+\sigma\,w_2(t)\end{array}\right)+{\mathcal A}_{\epsilon}(y(t))w(t)+{\mathcal B}_{\epsilon}(y(t))\left(\begin{array}{c}-w_2(t)\\w_1(t)\end{array}\right)\\
& = & \left( 
\begin{array}{cc}
\sigma & \mu\\
-\mu & \sigma
\end{array}
\right)	
\left(
\begin{array}{c}
w_1(t)\\
w_2(t)
\end{array}
\right)
)+{\mathcal A}_{\epsilon}(y(t))w(t)+{\mathcal B}_{\epsilon}(y(t))\left(\begin{array}{c}-w_2(t)\\w_1(t)\end{array}\right).
\end{eqnarray*}  

2.5. From $r(0)=\sqrt{(y_1(0))^2+(y_2(0))^2}$ and $\phi(0)=\psi$ we have 
$$
w(0)=r(0)\left(\begin{array}{c}\cos(\phi(0))\\\sin(\phi(0))\end{array}\right)=\sqrt{(y_1(0))^2+(y_2(0))^2}\left(\begin{array}{c}\cos(\psi)\\\sin(\psi)\end{array}\right)=\left(\begin{array}{c}y_1(0)\\y_2(0)\end{array}\right).
$$ 
So $w$ and $\left(
\begin{array}{c}
y_1\\
y_2
\end{array}
\right)$ both are solutions of the same initial value problem. Hence
$$
\left(
\begin{array}{c}
y_1(t)\\
y_2(t)
\end{array}
\right)=w(t)=r(t)\left(\begin{array}{c}\cos(\phi(t))\\\sin(\phi(t))\end{array}\right)\quad\mbox{for all}\quad t\in I.
$$

3. On assertion (iii). We have $(y_1(t))^2+(y_2(t))^2>0$
for all $t\in I$ since otherwise the invariance of $U$ under $V_{\epsilon}$ in $B_1$ would result in 
$y(t_0)\in U$, hence $y_1(t_0)=0=y_2(t_0)$, in contradiction to the hypothesis on $y$. $\Box$

\medskip

The maximal solutions of the system (2,11,12) with initial values
$$
y(0)=x\in\mathbb{R}^3,\quad r(0)\in(0,\infty),\quad\phi(0)\in\mathbb{R}
$$
define a continuously differentiable flow 
$$
G_{\epsilon}:\mathbb{R}\times\mathbb{R}^5\supset\,dom(G_{\epsilon})\to\mathbb{R}^5.
$$
Notice that for a solution $y:I\to\mathbb{R}^3$ of Eq. (2)  and for $\phi$ as in Proposition 3.4 (ii) we obtain
$$
\left(t,\left(\begin{array}{c}y(0)\\\sqrt{(y_1(0))^2+(y_2(0))^2}\\\phi(0)\end{array}\right)\right)\in\,dom(G_{\epsilon})\quad\mbox{and}\quad\phi(t)=G_{\epsilon,5}\left(t,\left(\begin{array}{c}y(0)\\\sqrt{(y_1(0))^2+(y_2(0))^2}\\\phi(0)\end{array}\right)\right)
$$
for all $t\in I$.

\begin{corollary}
Assume H$_1(\eta,\epsilon)$.
For every solution $y$ of Eq. (2) on an interval $I\ni0$ with  $y(I)\subset B_1$
and for $r,\phi$ as in Proposition 3.4 (ii) we have
$$
r'(t)\in\sigma\,r(t)+[-\eta,\eta]r(t)\quad\mbox{and}\quad
\phi'(t)\in-\mu+[-\eta,\eta]\quad\mbox{for all}\quad t\in I.
$$
\end{corollary}

{\bf Proof.} Let $t\in I$. From Corollary 3.3 we infer
\begin{equation*}
\left|\left(
\begin{array}{c}
R_{\epsilon,L,1}(y(t))\\
R_{\epsilon,L,2}(y(t))
\end{array}
\right)
\right|=|R_{\epsilon,L}(y(t))|\le\eta|P_Ly(t)|=\eta\sqrt{(y_1(t))^2+(y_2(t))^2}=\eta\,r(t).
\end{equation*}
Hence $|{\mathcal A}_{\epsilon}(y(t))|\le\eta$ and $|{\mathcal B}_{\epsilon}(y(t))|\le\eta$. Use Eqs. (11,12) to complete the proof. $\Box$

\section{A travel time from a subset of $M_I$ to $M_E$}

The next objective is to establish the travel time from the continuously differentiable submanifold 
$\{x\in M_I:0<x_3<1\}$ to the  differentiable submanifold $M_E$ as a continuously differentiable map. 
 
\begin{proposition}
Assume 
$$
\mbox{H}_2(\eta,\epsilon)\qquad 0<\eta<\min\{-\sigma/\sqrt{2},u/2\}\quad\mbox{and}\quad0<\epsilon<\min\{\epsilon(\eta),\epsilon_M,\epsilon_I,\epsilon_E,\epsilon(u/2),\epsilon(-\sigma/2)\}.
$$

(i) For every $x\in M_I$ with $0<x_3<1$ there exists $t>0$ with $(t,x)\in\Omega_{\epsilon}$ and $F_{\epsilon}(t,x)\notin B_1$.

\medskip

(ii) The map $t_{\epsilon}:\{x\in M_I:0<x_3<1\}\to\mathbb{R}$ given by
$$
t_{\epsilon}(x)=\inf\,\{t>0:(t,x)\in\Omega_{\epsilon}\,\,\mbox{and}\,\,F_{\epsilon}(t,x)\notin B_1\}
$$
is continuous. For every $x\in M_I$ with $0<x_3<1$ we have
$t_{\epsilon}(x)>0$, and the maximal solution $y$ of Eq. (2) with $y(0)=x$ satisfies $y_3(t_{\epsilon}(x))=1$, and $y_3'(t)>0$ for $0\le t\le t_{\epsilon}(x)$, and 
$0<x_3<y_3(t)<1$ for $0<t<t_{\epsilon}(x)$.

\medskip

(iii) The map $t_{\epsilon}$ is continuously differentiable.
\end{proposition}

{\bf Proof.} 1. On (i). For $x\in M_I$ with $0<x_3<1$ let $y:I_y\to\mathbb{R}^3$ denote the maximal solution of  Eq. (2) with  initial value $y(0)=x$, so that $y(t)=F(t,x)$ on $I_y$. Suppose $y(t)\in B_1$ for all $t\ge0$ in $I_y$. In case $\sup\, I_y=\infty$ we obtain by means of Proposition 3.4 (i) first
$y_3(t)>0$ for all $t\ge0$, and then $y_3(t)\ge x_3e^{(u-\eta)t}$ for all $t\ge0$, which contradicts the assumption of boundedness.  In case $\sup\, I_y<\infty$ we have $I_y\cap[0,\infty)=[0,\tau)$ for some $\tau>0$. By assumption $y(t)\in B_1$ on $[0,\tau)$. As $V_{\epsilon}$ is bounded on $B_1$, we get a bound for $y'(t)$ on $[0,\tau)$. Hence $y$ satisfies a Lipschitz condition on $[0,\tau)$. This yields that $y$ has a limit at $\tau$, which can be used to construct a continuation of $y$ as a solution of Eq. (2) beyond $\tau$, in contradiction to maximality.

\medskip

2. Let $x\in M_I$ with $0<x_3<1$ be given. Due to assertion (i) we obtain $t_{\epsilon}(x)\in[0,\infty)$. Consider $y:I_y\to\mathbb{R}^3$ as in Part 1. Proof  of 
$t_{\epsilon}(x)>0$ and of the assertions concerning $y$. We have $I_y\cap[0,\infty)=[0,\tau)$ with $0<\tau\le\infty$.
Obviously, $0\le t_{\epsilon}<\tau$. For $r=\bigg|\left(\begin{array}{c}y_1\\y_2\end{array}\right)\bigg|$
we have $r(0)=1$, and by Corollary 3.5, $r'(t)\le(\sigma+\eta)r(t)<0$ on $[0,t_{\epsilon}(x)]$. It follows that $r(t)<1$ on $(0,t_{\epsilon}(x)]$. Using $y(s)\notin B_1$ for some $s>t_{\epsilon}(x)$ arbitrarily close to $t_{\epsilon}(x)$ in combination with $r(s)<1$ we conclude that $|y_3(s)|>1$ for these arguments $s$. By continuity, $|y_3(t_{\epsilon}(x))|=1$. This yields $t_{\epsilon}(x)>0$, in view of $y_3(0)=x_3\in(0,1)$. Using $0<x_3=y_3(0)<1$ and Proposition 3.4 (i) we deduce that $y_3'(t)\ge(u-\eta)y_3(t)>0$ for $0\le t\le t_{\epsilon}(x)$, and upon that,
$0<y_3(t_{\epsilon}(x))$. It follows that $y_3(t_{\epsilon}(x))=1$. Also, $0<x_3<y_3(t)<1$ on $(0,t_{\epsilon}(x))$.

\medskip

3. On continuity of the map $t_{\epsilon}$. Let $x\in M_I$ with $0<x_3<1$ be given and consider $y$ and $r$ as in Part 2. Let $\rho\in(0,t_{\epsilon}(x))$ be given. 
For some $s\in(t_{\epsilon}(x),t_{\epsilon}(x)+\rho)$, $y(s)\notin B_1$. By continuity there is a neighbourhood $N_1$ of $x$ in $\mathbb{R}^3$ with $F_{\epsilon}(s,z)\notin B_1$ for all $z\in N_1$. According to Part 2, $0<y_3(t)<1$ and $0<r(t)<1$ on $[0,t_{\epsilon}(x)-\rho]$. Hence $F_{\epsilon}(t,x)\in\,int\,B_1$ on $[0,t_{\epsilon}(x)-\rho]$. By continuity and compactness we find a neighbourhood $N\subset N_1$ of $x$ in $\mathbb{R}^3$ so that for all $z\in N$  and for all $t\in[0,t_{\epsilon}(x)-\rho]$ we have  $F_{\epsilon}(t,z)\in\,int\,B_1$. Hence
$t_{\epsilon}(x)-\rho\le t_{\epsilon}(z)<s<\tau_{\epsilon}(x)+\rho$ for all $z\in N\subset N_1,$ which yields continuity of $t_{\epsilon}$ at $x$.

\medskip

4. We show that locally the map $t_{\epsilon}$ is given by continuously differentiable maps. Let $x\in M_I$ with $0<x_3<1$ be given. 
Then $F_{\epsilon,3}(t_{\epsilon}(x),x)=1$, and
$\partial_1F_{\epsilon,3}(t_{\epsilon}(x),x)>0$. The Implicit Function Theorem yields an open neighbourhood $N$ of $x$ in $\mathbb{R}^3$
and $\xi>0$ and a continuously differentiable map $\tau:N\to (t_{\epsilon}-\xi,t_{\epsilon}+\xi)$ 
with $F_{\epsilon,.3}(\tau(z),z)=1$ for all $z\in N$, and on
$(t_{\epsilon}-\xi,t_{\epsilon}+\xi)\times N$,
\begin{equation*}
F_{\epsilon,3}(t,z)=1\quad\mbox{if and only if}\quad t=\tau(z).
\end{equation*}
By continuity there is an open neighbourhood $N_1\subset N$ of $x$ in $\mathbb{R}^3$ so that for all $z\in N_1\cap M_I$ we have $0<z_3<1$ and $t_{\epsilon}(x)-\xi<t_{\epsilon}(z)<t_{\epsilon}(x)+\xi$.   Recall $F_{\epsilon,3}(t_{\epsilon}(z),z)=1$ for all $z\in M_I$ with $0<z_3<1$. It follows that on $N_1\cap M_I$, $t_{\epsilon}(z)=\tau(z)$.
The restriction of $\tau$ to $N_1\cap M_I$ is a continuously differentiable function on the open subset $N_1\cap M_I$ of the submanifold $M_I$. $\Box$

\begin{corollary}
Assume H$_2(\eta,\epsilon)$. For each $x\in M_I$ with $0<x_3<1$,
\begin{equation*}
\frac{1}{u+\eta}(-\log\,x_3)\le t_{\epsilon}(x)\le \frac{1}{u-\eta}(-\log\,x_3).
\end{equation*}
\end{corollary}

{\bf Proof.} Using Proposition 3.4  (i) and $x_3<F_{\epsilon,3}(t,x)<1=F_{\epsilon,3}(t_{\epsilon}(x),x)$ on $(0,t_{\epsilon}(x))$   we get
\begin{equation*}
1=F_{\epsilon,3}(t_{\epsilon}(x),x)\in x_3e^{(u+[-\eta,\eta])t_{\epsilon}(x)}
\end{equation*}
which yields the assertion. $\Box$

\section{Transversality at $M_I$ and the inner map}

The {\it inner map}
\begin{equation*}
I_{\epsilon}:\{z\in M_I:0<z_3<1\}\ni x\mapsto F_{\epsilon}(t_{\epsilon}(x),x)\in M_E,
\end{equation*}
for $0<\epsilon<\epsilon(\eta)$ and $0<\eta<\min\{-\frac{\sigma}{\sqrt{2}},\frac{u}{2}\}$,
is continuously differentiable, from the submanifold $M_I$ into the submanifold $M_E$. 

\begin{proposition}
Assume H$_2(\eta,\epsilon)$. The inner map defines a diffeomorphism onto the open subset $I_{\epsilon}(\{x\in M_I:0<x_3<1\})$ of $M_E$.
\end{proposition}

{\bf Sketch of the proof.} 1. The transversality condition
\begin{equation*}
D_1 F_{\epsilon}(t_{\epsilon}(x),x)1\notin L=T_{ F_{\epsilon}(t_{\epsilon}(x),x)}M_E
\end{equation*}
is satisfied since according to Proposition 4.1 (ii),
\begin{equation*}
[D_1 F_{\epsilon}(t_{\epsilon}(x),x)1]_3=y_3'(t_{\epsilon}(x))>0
\end{equation*}
for $y=F_{\epsilon}(\cdot,x)$ with $x\in M_I$ and $0<x_3<1$.

\medskip

2. Proof that for every $x\in M_I$ with $0<x_3<1$ the derivative $T_xI_{\epsilon}:T_xM_I\to T_{I_{\epsilon}(x)}M_E$ is an isomorphism. This derivative is given by
\begin{equation*}
T_xI_{\epsilon}v=PD_2F_{\epsilon}(t_{\epsilon}(x),x)v\quad\mbox{for all}\quad v\in T_xM_I,
\end{equation*}
with the projection $P:\mathbb{R}^3\to\mathbb{R}^3$ along $D_1 F_{\epsilon}(t_{\epsilon}(x),x)1$ onto $L=T_{ F_{\epsilon}(t_{\epsilon}(x),x)}M_E$. The isomorphism $D_2F_{\epsilon}(t_{\epsilon}(x),x)$ maps $T_xM_I$ onto a plane $Z$ and sends the
vector $D_1F_{\epsilon}(0,x)1=V_{\epsilon}(x)\notin T_xM_I$ (see Proposition 3.1 (ii)) to $D_1F_{\epsilon}(t_{\epsilon}(x),x)1$ which is transversal to $Z$. Therefore the projection $P$ defines an isomorphism from $Z$ onto $T_{F_{\epsilon}(t_{\epsilon}(x),x)}M_E$.

\medskip

3. The Inverse Mapping Theorem yields that the inner map is locally given by diffeomorphisms from open subsets of $M_I$ onto open subsets of $M_E$. 

\medskip

4. The assertion follows easily provided the inner map is injective. Proof of this: Suppose $ I_{\epsilon}(x)= I_{\epsilon}(z)$ for $x,z$ in $M_I$ with $0<x_3<1,0<z_3<1$.
Let $y=F_{\epsilon}(\cdot,x)$ and $w=F_{\epsilon}(\cdot,z)$, so that $y(t_{\epsilon}(x))=w(t_{\epsilon}(z))$. Suppose
$t_{\epsilon}(x)\neq t_{\epsilon}(z)$. In case $t_{\epsilon}(x)<t_{\epsilon}(z)$ uniqueness for backward solutions of initial value problems yields $y(t_{\epsilon}(x)-s)=w(t_{\epsilon}(z)-s)$ for all $s\in[0,t_{\epsilon}(x)]$. Hence $x=y(0)=y(t_{\epsilon}(x)-t_{\epsilon}(x))=w(t_{\epsilon}(z)-t_{\epsilon}(x))$, and thereby $1=|P_Lx|=|P_Lw(t_{\epsilon}(z)-t_{\epsilon}(x))|$, in contradiction to $|P_Lz|=1$ and the fact that due to Corollary 3.5  $|P_Lw|$ is strictly decreasing on $[0,t_{\epsilon}(z)]$. It follows that $t_{\epsilon}(x)=t_{\epsilon}(z)$, from which we infer $x=y(0)=y(t_{\epsilon}(x)-t_{\epsilon}(x))=w(t_{\epsilon}(z)-t_{\epsilon}(x))=w(0)=z$. $\Box$

\medskip

It will be convenient to describe the inner map (and the outer map, which will be introduced in the next section) in terms of the parametrization 
\begin{equation*}
K_{\epsilon}:(-\pi,\pi)\times\mathbb{R}\ni(\psi,\delta)\mapsto\left(\begin{array}{c}\cos(\omega_{\epsilon}+\psi) \\ \sin(\omega_{\epsilon}+\psi) \\ \delta\end{array}\right)\in\mathbb{R}^3,\quad 0<\epsilon<\epsilon_M,
\end{equation*}
where $\omega_{\epsilon}\in[0,2\pi)$ is chosen so that
\begin{equation*}
\left(\begin{array}{c}\cos(\omega_{\epsilon}) \\ \sin(\omega_{\epsilon}) \\ 0\end{array}\right)=h_{\epsilon}(t_{I,\epsilon})\in M_I\cap L.
\end{equation*}
 The map $K_{\epsilon}$ defines a diffeomorphism from $(-\pi,\pi)\times(0,1)$ onto the open subset 
\begin{equation*}
M_I\setminus\left\{\left(\begin{array}{c}\cos(\omega_{\epsilon}+\pi) \\ \sin(\omega_{\epsilon}+\pi) \\ \delta\end{array}\right)\in\mathbb{R}^3:\delta\in\mathbb{R}\right\}
\end{equation*}
of $M_I$. Compare Figure 3 on page 15.

\begin{figure}
	\includegraphics[page=1,scale=0.7]{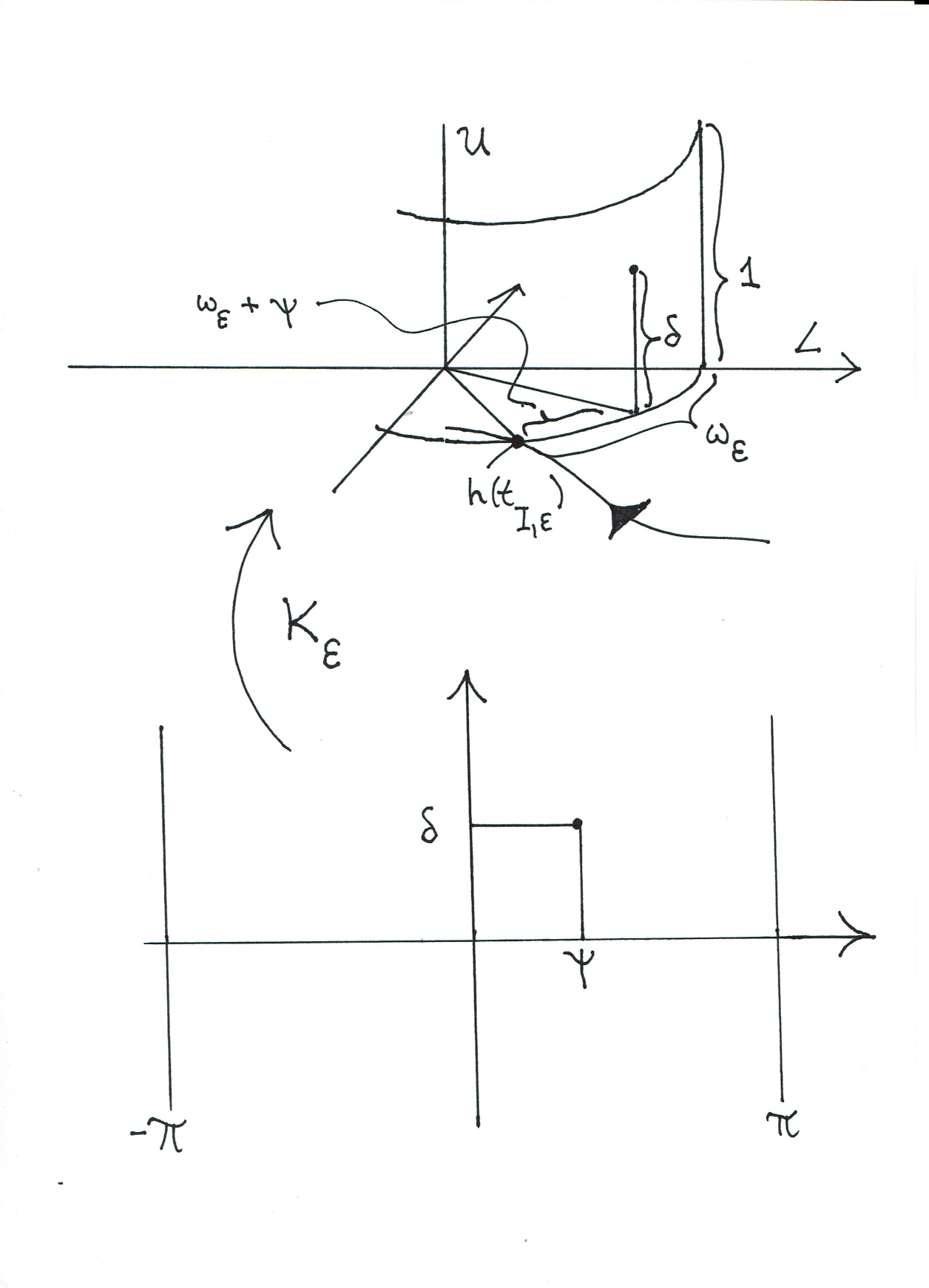}
	\caption{Parametrizing the domain of the inner map}
\end{figure}

We turn to the angle of the projection into the plane $L$ of the value of the inner map. Notice the difference to angles for projected solutions of Eq. (2), which were addressed in Section 3 ! For $x\in M_I$ with $0<x_3<1$ we have
$\sqrt{x_1^2+x_2^2}=1$, and a unique $\psi=\psi_x\in[\omega_{\epsilon}-\pi,\omega_{\epsilon}+\pi)$ with
$$
\left(\begin{array}{c}x_1\\x_2\end{array}\right)=\left(\begin{array}{c}\cos(\psi)\\\sin(\psi)\end{array}\right),
$$
and the maximal solution $y$ of Eq. (2) with $y(0)=x$ satisfies
$y(t)\in B_1$ on the interval $[0,t_{\epsilon}(x)]$, see Proposition 4.1 (ii). The function $\phi:[0,t_{\epsilon}(x)]\to\mathbb{R}$ given by $y$ and $\phi(0)=\psi=\psi_x$ as in Proposition 3.4 (ii) is uniquely determined by $x$, so we set $\phi_x=\phi$ and obtain the map
$$
\Phi_{\epsilon,\ast}:\{x\in M_I:0<x_3<1\}\ni x\mapsto \phi_x(t_{\epsilon}(x))\in\mathbb{R}.
$$
In view of the equation
$$
\left(\begin{array}{c}y_1(t)\\y_2(t)\end{array}\right)=\sqrt{(y_1(t))^2+(y_2(t))^2}\left(\begin{array}{c}\cos(\phi_x(t))\\\sin(\phi_x(t))\end{array}\right)\quad\mbox{with}\quad r_x(t)=\sqrt{(y_1(t))^2+(y_2(t))^2}
$$
at $t=t_{\epsilon}(x)$, by Proposition 3.4 (ii), the value $\Phi_{\epsilon,\ast}(x)$ 
stands for the desired angle. Compare Figure 4 on page 16.

\begin{figure}
	\includegraphics[page=1,scale=0.7]{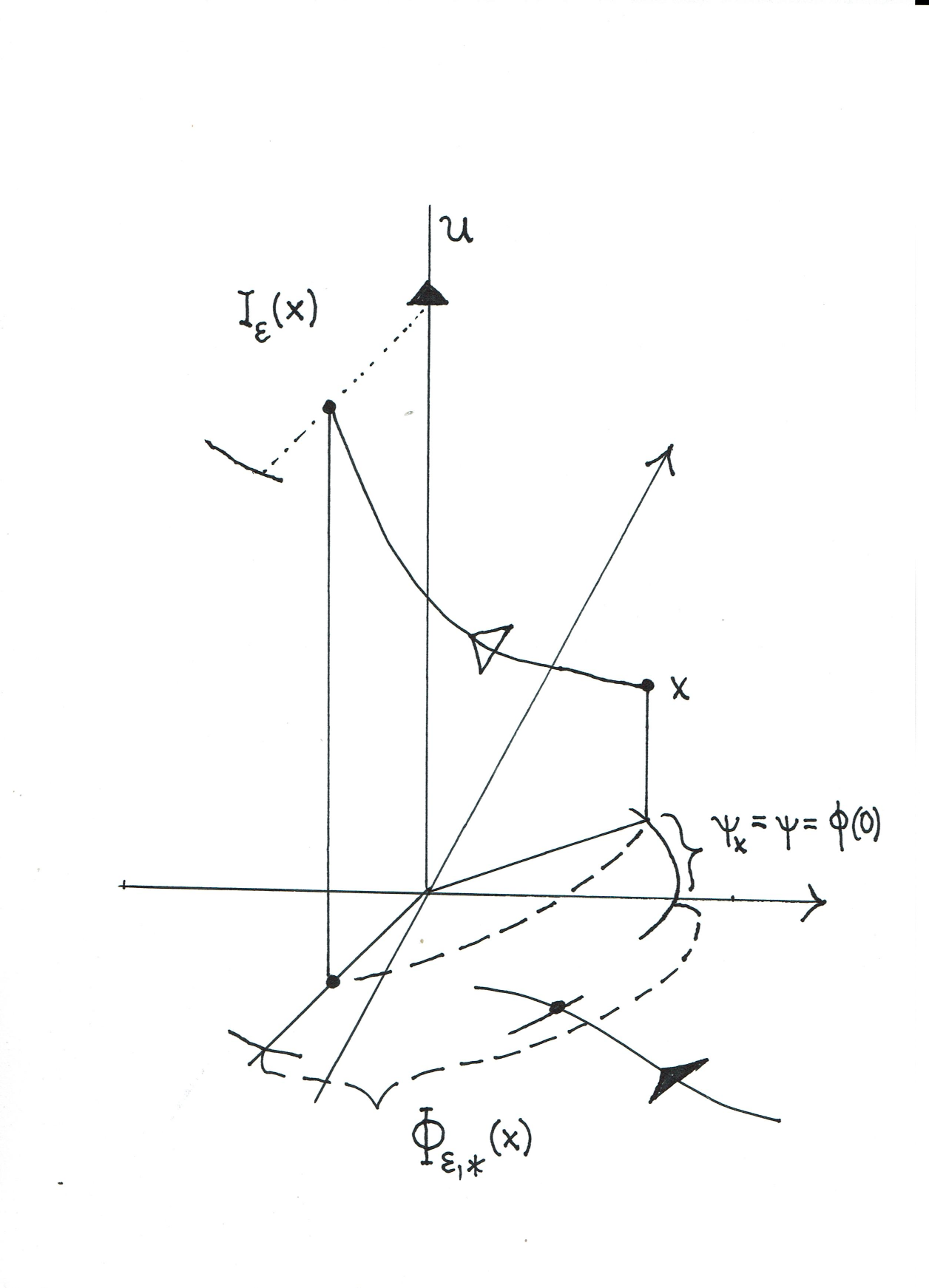}
	\caption{The action of the inner map and related angles for projections}
\end{figure}

\begin{proposition}
Assume H$_2(\eta,\epsilon)$.
The function $\Phi_{\epsilon}:(-\pi,\pi)\times(0,1)\to\mathbb{R}$ given by $\Phi_{\epsilon}(\psi,\delta)=\Phi_{\epsilon,\ast}(K_{\epsilon}(\psi,\delta))$ is continuously differentiable. 
\end{proposition}

{\bf Proof.} 1. By the remarks before Corollary 3.5,
$$
\Phi_{\epsilon,\ast}(x)=\phi_x(t_{\epsilon}(x))=G_{\epsilon,5}\left(t_{\epsilon}(x),\left(\begin{array}{c} x\\1\\\psi_x\end{array}\right)\right)
$$
for all $x\in M_I$ with $0<x_3<1$. The maps $t_{\epsilon}$ and $G_{\epsilon}$ are continuously differentiable, as well as the map
$$
\mathbb{R}^2\setminus\left\{r\left(\begin{array}{c}h_{\epsilon}(t_{I,\epsilon})_1\\h_{\epsilon}(t_{I,\epsilon})_2\end{array}\right):r\le0\right\}\ni\left(\begin{array}{c}x_1\\x_2\end{array}\right)\mapsto\psi_x\in(\omega_{\epsilon}-\pi,\omega_{\epsilon}+\pi)
$$
given by 
$$
\left(\begin{array}{c}x_1\\x_2\end{array}\right)=\sqrt{x_1^2+x_2^2}\left(\begin{array}{c}\cos(\psi_x)\\\sin(\psi_x)\end{array}\right).
$$
It follows that the real function $\Phi_{\epsilon,\ast}$ on the continuously differentiable submanifold $\{x\in M_I:0<x_3<1\,\,\mbox{and}\,\,P_Lx\neq-h_{\epsilon}(t_{I,\epsilon})\}$ is continuously differentiable. 

\medskip

2. The chain rule yields that also the map $\Phi_{\epsilon}$ is continuously differentiable. $\Box$

\medskip

Next we estimate the range of the inner map in terms of the parameters in $(-\pi,\pi)\times(0,1)$.

\begin{proposition}
Assume H$_2(\eta,\epsilon)$. Let  $x=K_{\epsilon}(\psi,\delta)$ with
$-\pi<\psi<\pi$ and $0<\delta<1$ be given. 
	
\medskip
	
(i) Then $x\in M_I$ and $0<x_3=\delta<1$, and
\begin{equation*}
I_{\epsilon}(x)\in\left(
\begin{array}{c}e^{(\sigma+[-\eta,\eta])\left[\frac{1}{u+\eta},\frac{1}{u-\eta}\right]\log\left(\frac{1}{\delta}\right)}\left(\begin{array}{c}\cos\\ \sin\end{array}\right)\left(\omega_{\epsilon}+\psi+(-\mu+[-\eta,\eta])\left[\frac{1}{u+\eta},\frac{1}{u-\eta}\right]\log\left(\frac{1}{\delta}\right)\right) \\ 1 \end{array}\right),
\end{equation*}
	
(ii) $\qquad\delta^{\frac{-\sigma+\eta}{u-\eta}}\le|P_LI_{\epsilon}(K_{\epsilon}(\psi,\delta))|\le\delta^{\frac{-\sigma-\eta}{u+\eta}}.$

\end{proposition}

{\bf Proof.} 1. On assertion (i). We have 
$$
\left(\begin{array}{c}I_{\epsilon,1}(x)\\I_{\epsilon,2}(x)\end{array}\right) =r_x(t_{\epsilon}(x))\left(\begin{array}{c}\cos(\Phi_{\epsilon}(\psi,\delta))\\ \sin(\Phi_{\epsilon}(\psi,\delta))\end{array}\right),
$$
see the remarks prior to Proposition 5.3. By Corollary 4.2, $t_{\epsilon}(x)\in\left[\frac{1}{u+\eta},\frac{1}{u-\eta}\right]\log\left(\frac{1}{\delta}\right)$. Corollary 3.5, with $r(0)=1$ and $\phi(0)=\omega_{\epsilon}+\psi$, yields
\begin{eqnarray*}
r_x(t_{\epsilon}(x)) & \in &  e^{(\sigma+[-\eta,\eta])t_{\epsilon}(x)},\\
\Phi_{\epsilon}(\psi,\delta) & \in & \omega_{\epsilon}+\psi+(-\mu+[-\eta,\eta])t_{\epsilon}(x).
\end{eqnarray*}
Combining these relations and $I_{\epsilon,3}(x)=1$ we obtain assertion (i). 

\medskip

2. On assertion (ii). Using
$$	|P_LI_{\epsilon}(K_{\epsilon}(\psi,\delta))|=\bigg|\left(\begin{array}{c}I_{\epsilon,1}(x)\\I_{\epsilon,2}(x)\end{array}\right)\bigg|
$$
we infer from assertion (i) the estimates
\begin{eqnarray*}
\delta^{\frac{-\sigma+\eta}{u-\eta}} & = & 
\left(\frac{1}{\delta}\right)^{\frac{\sigma-\eta}{u-\eta}}=
e^{(\sigma-\eta)\log\left(\frac{1}{\delta}\right)\frac{1}{u-\eta}} \le 	
|P_LI_{\epsilon}(K_{\epsilon}(\psi,\delta))|
\le  e^{(\sigma+\eta)\log\left(\frac{1}{\delta}\right)\frac{1}{u+\eta}}\\
& = & \left(\frac{1}{\delta}\right)^{\frac{\sigma+\eta}{u+\eta}}=\delta^{\frac{-\sigma-\eta}{u+\eta}}.\quad\Box
\end{eqnarray*}

\begin{figure}
	\includegraphics[page=1,scale=0.7]{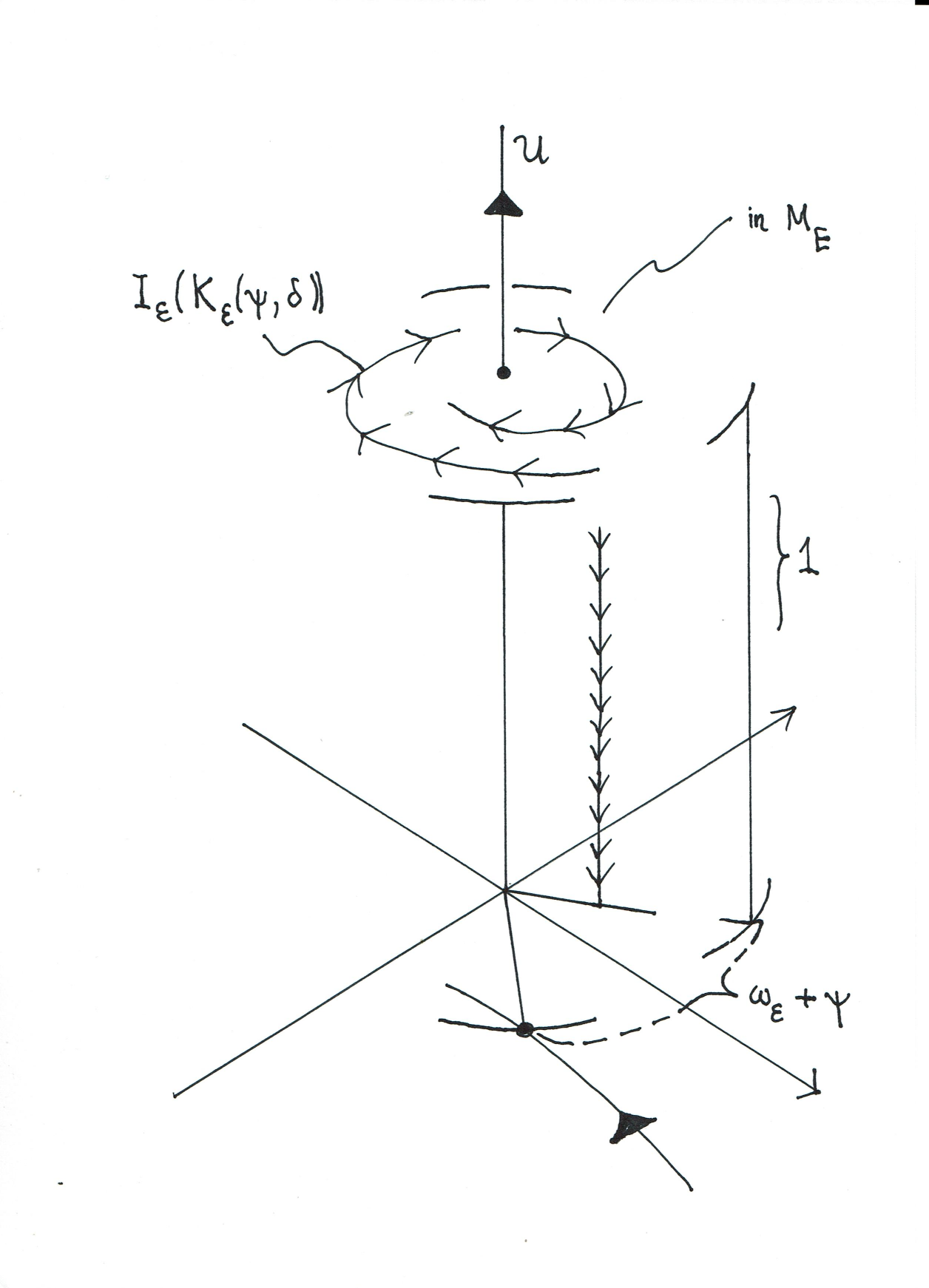}
	\caption{The inner map along vertical line segments}
\end{figure}

Notice that in vertical line segments, for $\psi$ fixed and $\delta\searrow0$, the angle $\Phi_{\epsilon}(\psi,\delta)$ decreases with nearly constant speed close to $-\mu<0$  to $-\infty$ and the radius $r_x(t_{\epsilon}(x))$, $x=K_{\epsilon}(\psi,\delta)$, decreases almost like $\delta^{-\sigma/u}$ to $0$, so $P_LI_{\epsilon}(x)$ spirals clockwise down to $0\in L$. 

\section{The outer map, and the return map in the plane}

Under condition  H$_2(\eta,\epsilon)$ the transversality property $V_{\epsilon}(x)\notin T_xM_I$ from Proposition 3.1 (ii), for $x=h_{\epsilon}(t_{I,\epsilon})\in L$ with $|x|=1$, yields a continuously differentiable travel time
$$
\tau_{\epsilon}:\{x\in M_E:|P_Lx|<r_{\epsilon}\}\to(0,\infty)
$$
for some $r_{\epsilon}>0$, with
$$
\tau_{\epsilon}(h_{\epsilon}(t_{E,\epsilon}))=t_{I,\epsilon}
\quad\mbox{and}\quad F_{\epsilon}(\tau_{\epsilon}(x),x)\in M_I
\quad\mbox{for all}\quad x\in M_E\quad\mbox{with}\quad|P_Lx|<r_{\epsilon}.
$$
The {\it outer map}
$$
E_{\epsilon}:\{x\in M_E:|P_Lx|<r_{\epsilon}\}\ni x\mapsto F_{\epsilon}(\tau_{\epsilon}(x),x)\in M_I
$$
is continuously differentiable, and we have
$$
E_{\epsilon}(h_{\epsilon}(t_{E,\epsilon}))=E_{\epsilon}(e_3)=h_{\epsilon}(t_{I,\epsilon})=\left(\begin{array}
{c}\cos(\omega_{\epsilon})\\ \sin(\omega_{\epsilon})\\0\end{array}\right),
$$
compare Figure 1 on page 3. Because of the transversality condition $V_{e_3}\notin T_{e_3}M_E$ from Proposition 3.1 (i) we may assume in addition that $E_{\epsilon}$ defines a diffeomorphism onto an open subset of $M_I$.

\medskip

Incidentally, $\tau_{\epsilon}$ and $E_{\epsilon}$ are given by restrictions of continuously differentiable maps on an open subset of $\mathbb{R}^3$ to an open subset of the submanifold $M_E$.

\medskip

By continuity there are $\delta_{I,\epsilon}\in(0,1)$ and $\omega_{I,\epsilon}\in(0,\pi)$ so that the set
$$
\left\{\left(\begin{array}
{c}\cos(\phi)\\ \sin(\phi)\\ \delta\end{array}\right)\in\mathbb{R}^3:|\phi-\omega_{\epsilon}|<\omega_{I,\epsilon},\,\,|\delta|<\delta_{I,\epsilon}\right\}
$$
is contained in the image of $E_{\epsilon}$. 

\medskip

The next result guarantees that the outer map sends a small disk centered at $e_3$ into the domain of the inverse $K_{\epsilon}^{-1}$, and that the inner map composed with $K_{\epsilon}$ sends a narrow strip into that small disk.

\begin{proposition}
Assume H$_2(\eta,\epsilon)$.

\medskip 

(i) There exists $r_{\epsilon,1}\in(0,r_{\epsilon})$ with $E_{\epsilon}(y)\in K_{\epsilon}((-\omega_{I,\epsilon},\omega_{I,\epsilon})\times(-\delta_{I,\epsilon},\delta_{I,\epsilon}))$
for all $y\in M_E$ with $|P_Ly|<r_{\epsilon,1}$.

\medskip

(ii) There exists $\delta_{I,\epsilon,1}\in(0,\delta_{I,\epsilon})$ so that for all $\psi\in(-\omega_{I,\epsilon},\omega_{I,\epsilon})$ and for all $\delta\in(0,\delta_{I,\epsilon,1})$ we have 
$$
\delta^{\frac{-\sigma+\eta}{u-\eta}}\le
|P_LI_{\epsilon}(K_{\epsilon}(\psi,\delta))|\le\delta^{\frac{-\sigma-\eta}{u+\eta}}<r_{\epsilon,1}.
$$
\end{proposition}

{\bf Proof.} 1. The set $K_{\epsilon}((-\omega_{I,\epsilon},\omega_{I,\epsilon})\times(-\delta_{I,\epsilon},\delta_{I,\epsilon}))$
is an open neighbourhood of $K_{\epsilon}(0,0)=h_{\epsilon}(t_{I,\epsilon})=E_{\epsilon}(e_3)$ in $M_I$. By continuity $E_{\epsilon}$ maps an open neighbourhood of $e_3$ in $M_E$ into $K_{\epsilon}((-\omega_{I,\epsilon},\omega_{I,\epsilon})\times(-\delta_{I,\epsilon},\delta_{I,\epsilon}))$.
The said neighbourhood of $e_3$ contains sets of the form $\{y\in M_E:|P_Ly|<r_{\epsilon,1}\}$ for $r_{\epsilon,1}\in (0,r_{\epsilon})$ sufficiently small.

\medskip

2. On (ii). Recall  Proposition 5.3 (ii) and choose $\delta_{I,\epsilon,1}\in(0,\delta_{I,\epsilon})$ so small that $\delta^{\frac{-\sigma-\eta}{u+\eta}}<r_{\epsilon,1}$ for $0<\delta<\delta_{I,\epsilon,1}$
which is possible due to  $0<\frac{-\sigma-\eta}{u+\eta}$. $\Box$

\medskip

Proposition 6.1 shows that the expression
$$
K_{\epsilon}^{-1}(E_{\epsilon}(I_{\epsilon}(K_{\epsilon}(\psi,\delta))))
$$
defines a diffeomorphism $R_{\epsilon}$ from $(-\omega_{I,\epsilon},\omega_{I,\epsilon})\times(0,\delta_{I,\epsilon,1})\subset\mathbb{R}^2$
onto an open subset of $(-\omega_{I,\epsilon},\omega_{I,\epsilon})\times(-\delta_{I,\epsilon},\delta_{I,\epsilon})\subset\mathbb{R}^2$.

\medskip

In addition to the previous representation of the return map 
$x\mapsto E_{\epsilon}(I_{\epsilon}(x))$ as a diffeomorphism $R_{\epsilon}$ in the plane we need a coordinate representation of the outer map alone.

\medskip

The expression $K_{\epsilon}^{-1}(E_{\epsilon}(x))$ defines a diffeomorphism $K_{\epsilon}^{-1}(E_{\epsilon}(\cdot))$ from the open subset $\{x\in M_E:|P_Lx|<r_{\epsilon,1}\}$ of $M_E$ onto an open subset of $\mathbb{R}^2$ which contains $(0,0)=K_{\epsilon}^{-1}(E_{\epsilon}(e_3))$. The linear map
$D(K_{\epsilon}^{-1}(E_{\epsilon}(\cdot)))(e_3)$ is an isomorphism from $T_{e_3}M_E=L$ onto $\mathbb{R}^2$. Let $v_{\epsilon}$ and $w_{\epsilon}$ denote the preimages under this isomorphism of the unit vectors $\left(\begin{array}{c} 1 \\ 0\end{array}\right)\in\mathbb{R}^2$ and $\left(\begin{array}{c} 0 \\ 1\end{array}\right)\in\mathbb{R}^2$, respectively.

\medskip

Let $\kappa_{\epsilon}:L\to\mathbb{R}^2$ denote the isomorphism given by $\kappa v_{\epsilon}=\left(\begin{array}{c} 1 \\ 0\end{array}\right)\in\mathbb{R}^2$ and $\kappa_{\epsilon} w_{\epsilon}=\left(\begin{array}{c} 0 \\ 1\end{array}\right)\in\mathbb{R}^2$.

\medskip

The restriction $P_{LE}=P_L|M_E$ which only subtracts $e_3$ is 
a diffeomorphism onto the plane with $DP_{LE}(y)\hat{y}=\hat{y}$ for all $y\in M_E$ and $\hat{y}\in T_yM_E=L$.

\medskip

The composition $\kappa_{\epsilon}\circ P_{LE}$ is a diffeomorphism from 
$M_E$ onto the plane, with  $\kappa_{\epsilon} P_{LE}(e_3)=\kappa_{\epsilon}(0,0)=(0,0)$.

\medskip

By continuity there exists $r_{E,\epsilon}>0$ so that $(\kappa_{\epsilon}\circ P_{LE})^{-1}=P_{LE}^{-1}\circ\kappa_{\epsilon}^{-1}$ maps
$\{(\psi,\delta)\in\mathbb{R}^2:|(\psi,\delta)|< r_{E,\epsilon}\}$ into the set $\{y\in M_E:|P_Ly|<r_{\epsilon,1}\}$. It follows that the expression
$$
K_{\epsilon}^{-1}(E_{\epsilon}((\kappa_{\epsilon}\circ P_{LE})^{-1}(\psi,\delta)))
$$
defines a diffeomorphism $K_{\epsilon}^{-1}(E_{\epsilon}((\kappa_{\epsilon}\circ P_{LE})^{-1})(\cdot,\cdot)))$ from $\{(\psi,\delta)\in\mathbb{R}^2:|(\psi,\delta)|< r_{E,\epsilon}\}$
onto an open subset of $\mathbb{R}^2$, with fixed point $(0,0)$.

\begin{corollary}
Assume H$_2(\eta,\epsilon)$. We have
$$
D(K_{\epsilon}^{-1}(E_{\epsilon}((\kappa_{\epsilon}\circ P_{LE})^{-1})(\cdot,\cdot))))(0,0)=id_{\mathbb{R}^2}.
$$ 
\end{corollary}

{\bf Proof.} Use the chain rule in combination with the previous statements about $\kappa_{\epsilon}$ and about the derivatives of $K_{\epsilon}^{-1}(E_{\epsilon}(\cdot))$ and $P_{LE}$. $\Box$

\medskip

The next result makes precise in which sense the diffeomorphism $K_{\epsilon}^{-1}(E_{\epsilon}((\kappa_{\epsilon}\circ P_{LE})^{-1})(\cdot,\cdot)))$ is close to the identity map.

\begin{corollary}
Assume H$_2(\eta,\epsilon)$.
For every $\beta>0$ there exists $\alpha_{\beta,\epsilon}>0$ so that for all 
$(\psi,\delta)\in[-\alpha_{\beta,\epsilon},\alpha_{\beta,\epsilon}]\times[-\alpha_{\beta,\epsilon},\alpha_{\beta,\epsilon}]$
we have	
\begin{eqnarray}
|(\psi,\delta)| & < & r_{E,\epsilon},\nonumber\\
|K_{\epsilon}^{-1}(E_{\epsilon}((\kappa_{\epsilon}\circ P_{LE})^{-1})(\psi,\delta)))-(\psi,\delta)| & \le & \beta|(\psi,\delta)|,\\
|D(K_{\epsilon}^{-1}(E_{\epsilon}((\kappa_{\epsilon}\circ P_{LE})^{-1}))(\psi,\delta)-id_{\mathbb{R}^2}|
& \le & \beta.\nonumber
\end{eqnarray}
\end{corollary}

{\bf Proof.} Use the definition of differentiability for $K_{\epsilon}^{-1}(E_{\epsilon}((\kappa_{\epsilon}\circ P_{LE})^{-1})(\cdot,\cdot))$ at $(0,0)$ and Corollary 6.2, and continuity of the derivative at $(0,0)$. $\Box$

\medskip

We turn to estimates of the range of the maps $\kappa_{\epsilon} P_{LE}(I_{\epsilon}(K_{\epsilon}(\cdot,\cdot)))$ and $R_{\epsilon}$.
The hypothesis H$_3(\eta,\epsilon)$ in the following proposition is stronger than H$_2(\eta,\epsilon)$ since $\eta<-\frac{\sigma}{2}$ implies $\eta<-\frac{\sigma}{\sqrt{2}}$ and $\eta<\frac{u}{2}$.

\begin{proposition}
Assume 
$$
\mbox{H}_3(\eta,\epsilon)\qquad 0<\eta<\min\{\mu,-\sigma/2,u/2\}\quad\mbox{and}\quad0<\epsilon<\min\{\epsilon(\eta),\epsilon_M,\epsilon_I,\epsilon_E,\epsilon(u/2),\epsilon(-\sigma/2)\}.
$$
Let $0<\beta\le\frac{1}{2}$. There exists $\alpha_{\beta,\epsilon}>0$ with $\alpha_{\beta,\epsilon}<\min\{\omega_{I,\epsilon},\delta_{I,\epsilon,1}\}$ 
such that for
$$
\delta_{\beta,\epsilon}=\left(\frac{2}{3(|\kappa_{\epsilon}|+1)}\alpha_{\beta,\epsilon}\right)^{\frac{3u}{-\sigma}}
$$
we have $\delta_{\beta,\epsilon}\le\frac{2}{3}\alpha_{\beta,\epsilon}$
and for all $(\psi,\delta)\in[-\alpha_{\beta,\epsilon},\alpha_{\beta,\epsilon}]\times(0.\delta_{\beta,\epsilon}]$,
\begin{eqnarray}
|\kappa_{\epsilon}P_{LE}(I_{\epsilon}(K_{\epsilon}(\psi,\delta)))| & \le \frac{2}{3}\alpha_{\beta,\epsilon},\\
R_{\epsilon}(\psi,\delta) & \in & [-\alpha_{\beta,\epsilon},\alpha_{\beta,\epsilon}]\times[-\alpha_{\beta,\epsilon},\alpha_{\beta,\epsilon}].
\end{eqnarray}
\end{proposition}

{\bf Proof.} 1. In order to show the inequality $\delta_{\beta,\epsilon}\le\frac{2}{3}\alpha_{\beta,\epsilon}$ observe that $0<\eta<\min\{-\frac{\sigma}{2},\frac{u}{2}\}$ yields 
$$
\frac{u+\eta}{-\sigma-\eta}<\frac{u+\frac{u}{2}}{-\sigma+\frac{\sigma}{2}}=\frac{\frac{3u}{2}}{\frac{-\sigma}{2}}=\frac{3u}{-\sigma}
$$
It follows that
$$
\delta_{\beta,\epsilon}=\left(\frac{2}{3(|\kappa_{\epsilon}|+1)}\alpha_{\beta,\epsilon}\right)^{\frac{3u}{-\sigma}}\le\left(\frac{2}{3(|\kappa_{\epsilon}|+1)}\alpha_{\beta,\epsilon}\right)^\frac{u+\eta}{-\sigma-\eta}.
$$
Consequently, with $0<-\sigma-\eta<u+\eta$ and $\delta_{\beta,\epsilon}<1$,
$$
\delta_{\beta,\epsilon}\le(|\kappa_{\epsilon}|+1)\delta_{\beta,\epsilon}^{\frac{-\sigma-\eta}{u+\eta}}\le\frac{2}{3}\alpha_{\beta,\epsilon}.
$$

2. From Corollary 6.3  we get $\alpha_{\beta,\epsilon}>0$ with $\alpha_{\beta,\epsilon}<\min\{\omega_{I,\epsilon},\delta_{I,\epsilon,1}\}$
so that for all $(\tilde{\psi},\tilde{\delta})\in\mathbb{R}^2$ with 
$(\tilde{\psi},\tilde{\delta})\in[-\alpha_{\beta,\epsilon},\alpha_{\beta,\epsilon}]\times[-\alpha_{\beta,\epsilon},\alpha_{\beta,\epsilon}]$
we have	
$|(\tilde{\psi},\tilde{\delta})|<r_{E,\epsilon}$ and, by means of (13) and with $0<\beta\le\frac{1}{2}$,
$$
|K_{\epsilon}^{-1}(E_{\epsilon}((\kappa_{\epsilon}\circ P_{LE})^{-1})(\tilde{\psi},\tilde{\delta})))|\le(1+\beta)|(\tilde{\psi},\tilde{\delta})|\le\frac{3}{2}|(\tilde{\psi},\tilde{\delta})|.
$$

3. Proof of (14). Because of $\delta_{\beta,\epsilon}\le\alpha_{\beta,\epsilon}<\delta_{I,\epsilon,1}$ we obtain from Proposition 6.1(ii) in combination with the result of Part 1 that for all $(\psi,\delta)\in(-\omega_{I,\epsilon},\omega_{I,\epsilon})\times(0,\delta_{\beta,\epsilon}]$, 
$$
|\kappa_{\epsilon} P_{LE}(I_{\epsilon}(K_{\epsilon}(\psi,\delta)))|\le
|\kappa_{\epsilon}|\delta^{\frac{-\sigma-\eta}{u+\eta}}\le
|\kappa_{\epsilon}|\delta_{\beta,\epsilon}^{\frac{-\sigma-\eta}{u+\eta}}\le\frac{2}{3}\alpha_{\beta,\epsilon}.
$$ 

4. Proof of (15). For $(\psi,\delta)$ as  in Part 3 let $(\tilde{\psi},\tilde{\delta})=\kappa P_{LE}(I_{\epsilon}(K_{\epsilon}(\psi,\delta)))$. Then
$$
|(\tilde{\psi},\tilde{\delta})|\le\frac{2}{3}\alpha_{\beta,\epsilon},
$$
hence $(\tilde{\psi},\tilde{\delta})\in[-\alpha_{\beta,\epsilon},\alpha_{\beta,\epsilon}]\times[-\alpha_{\beta,\epsilon},\alpha_{\beta,\epsilon}]$, which according to Part 2 yields $|(\tilde{\psi},\tilde{\delta})|<r_{E,\epsilon}$ and 
$$
|K^{-1}_{\epsilon}(E_{\epsilon}((\kappa_{\epsilon}\circ  P_{LE})^{-1}(\tilde{\psi},\tilde{\delta})))|\le \frac{3}{2}|(\tilde{\psi},\tilde{\delta})|.
$$
It follows that
\begin{eqnarray*}
|K^{-1}_{\epsilon}(E_{\epsilon}(I_{\epsilon}(K_{\epsilon}(\psi,\delta))))|
& = & |K^{-1}_{\epsilon}(E_{\epsilon}((\kappa_{\epsilon}\circ  P_{LE})^{-1}[(\kappa_{\epsilon}\circ P_{LE})(I_{\epsilon}(K(\psi,\delta)))])|\\
& = & |K^{-1}_{\epsilon}(E_{\epsilon}((\kappa_{\epsilon}\circ  P_{LE})^{-1}(\tilde{\psi},\tilde{\delta})))|\\
& \le & \frac{3}{2}|(\tilde{\psi},\tilde{\delta})|\le\alpha_{\beta,\epsilon}.
\end{eqnarray*}
Finally, use that the circle of radius $\alpha_{\beta,\epsilon}$ centered at the origin is contained in the square $[-\alpha_{\beta,\epsilon},\alpha_{\beta,\epsilon}]\times[-\alpha_{\beta,\epsilon},\alpha_{\beta,\epsilon}]$. $\Box$ 

\section{Curves expanded by the return map}

We show that the return map in coordinates $R_{\epsilon}$ expands each curve connecting two levels in the domain of 
$R_{\epsilon}$, in such a way that the resulting curve intersects the domain at least in two disjoint sets. These sets will be related to - but not given by - positions in the left or right halfplanes. 

\medskip

We begin with the angles $\Phi_{\epsilon}(\psi,\delta)=\Phi_{\epsilon,\ast}(x)$ of the projection into the plane $L$ of the values $I_{\epsilon}(x)$, for $x=K_{\epsilon}(\psi,\delta)$ with $-\pi<\psi<\pi$ and $0<\delta<1$. 

\medskip

For $0<\eta<\min\{\mu,u\}$ we define
\begin{eqnarray*}
	c & = & c_{\eta}=\frac{(u+\eta)(\mu+\eta)}{(u-\eta)(\mu-\eta)}\\
	& > & 1,\\
	k & = & k_{\eta}=e^{-6\pi\frac{u+\eta}{\mu-\eta}}\\
	& < & 1,
\end{eqnarray*}
and for $\delta_2\in(0,1)$ we define
\begin{eqnarray*}
\delta_1 & = & \delta_{1,\delta_2,\eta}=k_{\eta}\delta_2^{c_{\eta}}\\
& < & \delta_2.
\end{eqnarray*}

\begin{proposition}
Assume H$_2(\eta,\epsilon)$. Let $0<\delta_2<1$ and $\delta_1=\delta_{1,\delta_2,\eta}$. Then
$$
4\,\pi\le\Phi_{\epsilon}(\psi,\delta_2)-\Phi_{\epsilon}(\tilde{\psi},\delta_1)
$$
for all $\psi,\tilde{\psi}$ in $(-\pi,\pi)$,
\end{proposition}

{\bf Proof.} Assume $-\pi<\psi<\pi$, $-\pi<\tilde{\psi}<\pi$. Using the estimate of the speed of angles in Corollary 3.5 and the estimate of intersection times in Corollary 4.2 we have
$$
\Phi_{\epsilon}(\psi,\delta_2)-(\omega_{\epsilon}+\psi)\ge
(-\mu-\eta)\cdot\frac{1}{u-\eta}\log\left(\frac{1}{\delta_2}\right).
$$
Using this and the upper estimate
$$
\Phi_{\epsilon}(\tilde{\psi},\delta_1)-(\omega_{\epsilon}+\tilde{\psi})\le(-\mu+\eta)\cdot\frac{1}{u+\eta}\log\left(\frac{1}{\delta_1}
\right)
$$
we obtain
\begin{eqnarray*}
\Phi_{\epsilon}(\psi,\delta_2)-\Phi_{\epsilon}(\tilde{\psi},\delta_1) & \ge & -2\pi+\frac{\mu+\eta}{u-\eta}\log(\delta_2)
	-\frac{\mu-\eta}{u+\eta}\log(\delta_1)\\
& 	\ge &-2\pi+\frac{\mu+\eta}{u-\eta}\log(\delta_2)
-\frac{\mu-\eta}{u+\eta}[\log(k)+c\log(\delta_2)]\\
& = & -2\pi+6\pi+\left(\frac{\mu+\eta}{u-\eta}-c\frac{\mu-\eta}{u+\eta}\right)\log(\delta_2)=4\pi.\quad\Box
\end{eqnarray*}

It follows that
$$
m_1=\max_{|\psi|\le\alpha_{\beta}}\Phi_{\epsilon}(\psi,\delta_1)\quad\mbox{and}\quad m_2=\min_{|\psi|\le\alpha_{\beta}}\Phi_{\epsilon}(\psi,\delta_2)
$$ 
satisfy $m_1+4\pi\le m_2$. There exists $\psi_{\epsilon}\in[m_1+\pi,m_2-\pi]$ with
$$
\left(\begin{array}{c}\cos(\psi_{\epsilon})\\ \sin(\psi_{\epsilon})\\0\end{array}\right)=\frac{1}{|w_{\epsilon}|}w_{\epsilon}.
$$

\begin{proposition}
(Angles along curves connecting vertical levels) Assume H$_2(\eta,\epsilon)$. Let a (continuous) curve $c:[a,b]\to(-\pi,\pi)\times(0,1)$ be given with $c(b)_2=\delta_2$ and $c(a)_2=\delta_{1,\delta_2,\eta}=\delta_1$. Then there exist  
$a'_0<b'_0\le a'_1<b'_1$ in $[a,b]$
such that
$$
\Phi_{\epsilon}(c(t))\in(\psi_{\epsilon}-\pi,\psi_{\epsilon})\quad\mbox{on}\quad(a'_0,b'_0),\quad\Phi_{\epsilon}(c(a'_0))=
\psi_{\epsilon}-\pi,\quad\Phi_{\epsilon}(c(b'_0))=\psi_{\epsilon},
$$
$$
\Phi_{\epsilon}(c(t))\in(\psi_{\epsilon},\psi_{\epsilon}+\pi)\quad\mbox{on}\quad(a'_1,b'_1),\quad\Phi_{\epsilon}(c(a'_1))=\psi_{\epsilon},\quad\Phi_{\epsilon}(c(b'_1))=\psi_{\epsilon}+\pi.
$$
\end{proposition}

Compare Figure 6 on page 24.

\medskip

\begin{figure}
	\includegraphics[page=1,scale=0.7]{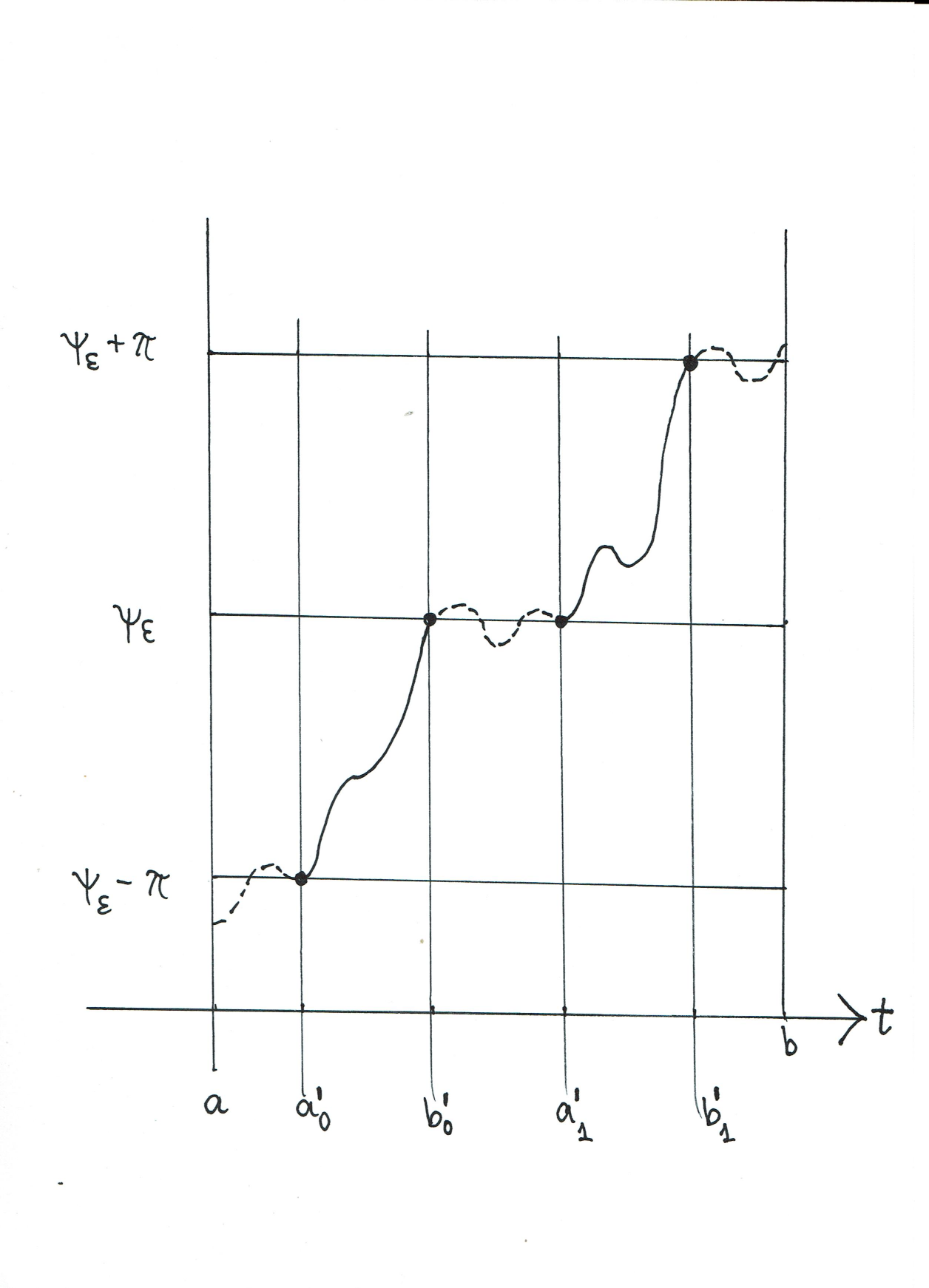}
	\caption{The map $\Phi_{\epsilon}\circ c$}
\end{figure}

{\bf Proof.} 1. We construct $a'_1$ and $b'_1$. From
$$
\Phi_{\epsilon}(c(a))\le m_1<m_1+\pi\le\psi_{\epsilon}\le m_2-\pi<m_2\le\Phi_{\epsilon}(c(b))
$$
we have
$$
\Phi_{\epsilon}(c(a))\le\psi_{\epsilon}-\pi<\psi_{\epsilon}<\psi_{\epsilon}+\pi\le\Phi_{\epsilon}(c(b)).
$$
By continuity, $\psi_{\epsilon}=\Phi_{\epsilon}(c(t))$ for some $t\in(a,b)$. Again by continuity there exists $b'_1\in(t,b]$ with 
$\Phi_{\epsilon}(c(s))<\psi_{\epsilon}+\pi$ on $[t,b'_1)$ and 
$\Phi_{\epsilon}(c(b'_1))=\psi_{\epsilon}+\pi$. Upon that, there exists
$a'_1\in[t,b'_1)$ with $\psi_{\epsilon}<\Phi_{\epsilon}(c(s))$ on $(a'_1,b'_1]$ and $\Phi_{\epsilon}(c(a'_1))=\psi_{\epsilon}$.

\medskip

2. The construction of $a'_0$ and $b'_0$ with $b'_0\le a'_1$ is analogous. $\Box$

\medskip

We turn to the height, or, to the second coordinate $R_{\epsilon,2}(\psi,\delta)$ of image points under the return map in coordinates, for arguments $(\psi,\delta)$ which under the inner map in coordinates 
$\kappa_{\epsilon} P_{LE}I_{\epsilon}(K_{\epsilon}(\cdot,\cdot))$
are mapped onto the vertical axis. Notice that in the cases
$$
\Phi_{\epsilon}(\psi,\delta)=\psi_{\epsilon}-\pi,\quad\Phi_{\epsilon}(\psi,\delta)=\psi_{\epsilon},\quad\Phi_{\epsilon}(\psi,\delta)=\psi_{\epsilon}+\pi
$$
we get that $P_{LE}I_{\epsilon}(K_{\epsilon}(\psi,\delta))$ belongs to the rays
$$
(0,\infty)(-w_{\epsilon}),\quad(0,\infty)w_{\epsilon},\quad(0,\infty)(-w_{\epsilon}),
$$ 
and $\kappa_{\epsilon}P_{LE}I_{\epsilon}(K_{\epsilon}(\psi,\delta))$ is on the vertical axis. 

\begin{proposition}
(From angles to vertical levels) Assume  H$_3(\eta,\epsilon)$ and 
\begin{equation}
0<\beta\le\frac{1}{2}.
\end{equation}
Consider 
\begin{equation}
\alpha_{\beta,\epsilon}<\min\{\omega_{I,\epsilon},\delta_{I,\epsilon,1}\}
\end{equation}
and  $\delta_{\beta,\epsilon}<\alpha_{\beta,\epsilon}$ according to Corollary 6.3 and Proposition 6.4. Assume
$\eta>0$ is so small that
\begin{equation}
c_{\eta}\frac{-\sigma+\eta}{u-\eta}<1,
\end{equation} 
and let  $\delta_2>0$ be so small that
\begin{equation}
\delta_2<\delta_{\beta,\epsilon}\quad\mbox{and}\quad
2\sqrt{2}\delta_2<\frac{1}{|\kappa_{\epsilon}^{-1}|}k_{\eta}^{\frac{-\sigma+\eta}{u-\eta}}\delta_2^{c_{\eta}\frac{-\sigma+\eta}{u-\eta}}.
\end{equation}
Let $\delta_1=\delta_{1,\delta_2,\eta}$ and let $(\psi,\delta)\in[-\alpha_{\beta,\epsilon},\alpha_{\beta,\epsilon}]\times[\delta_1,\delta_2]$ be given, and set $z=R_{\epsilon}(\psi,\delta)\in\mathbb{R}^2$.

\medskip

(i) Then $|z|>\sqrt{2}\delta_2$.

\medskip

(ii) In the cases 
$$
\Phi_{\epsilon}(\psi,\delta)=\psi_{\epsilon}-\pi,\quad\Phi_{\epsilon}(\psi,\delta)=\psi_{\epsilon},\quad \Phi_{\epsilon}(\psi,\delta)=\psi_{\epsilon}+\pi,
$$  
we have that
$$
z_2<-\delta_2,\quad z_2>\delta_2,\quad
z_2<-\delta_2,\quad\mbox{respectively}.
$$
\end{proposition}

Compare Figure 7 on page 26.

\medskip

\begin{figure}
	\includegraphics[page=1,scale=0.7]{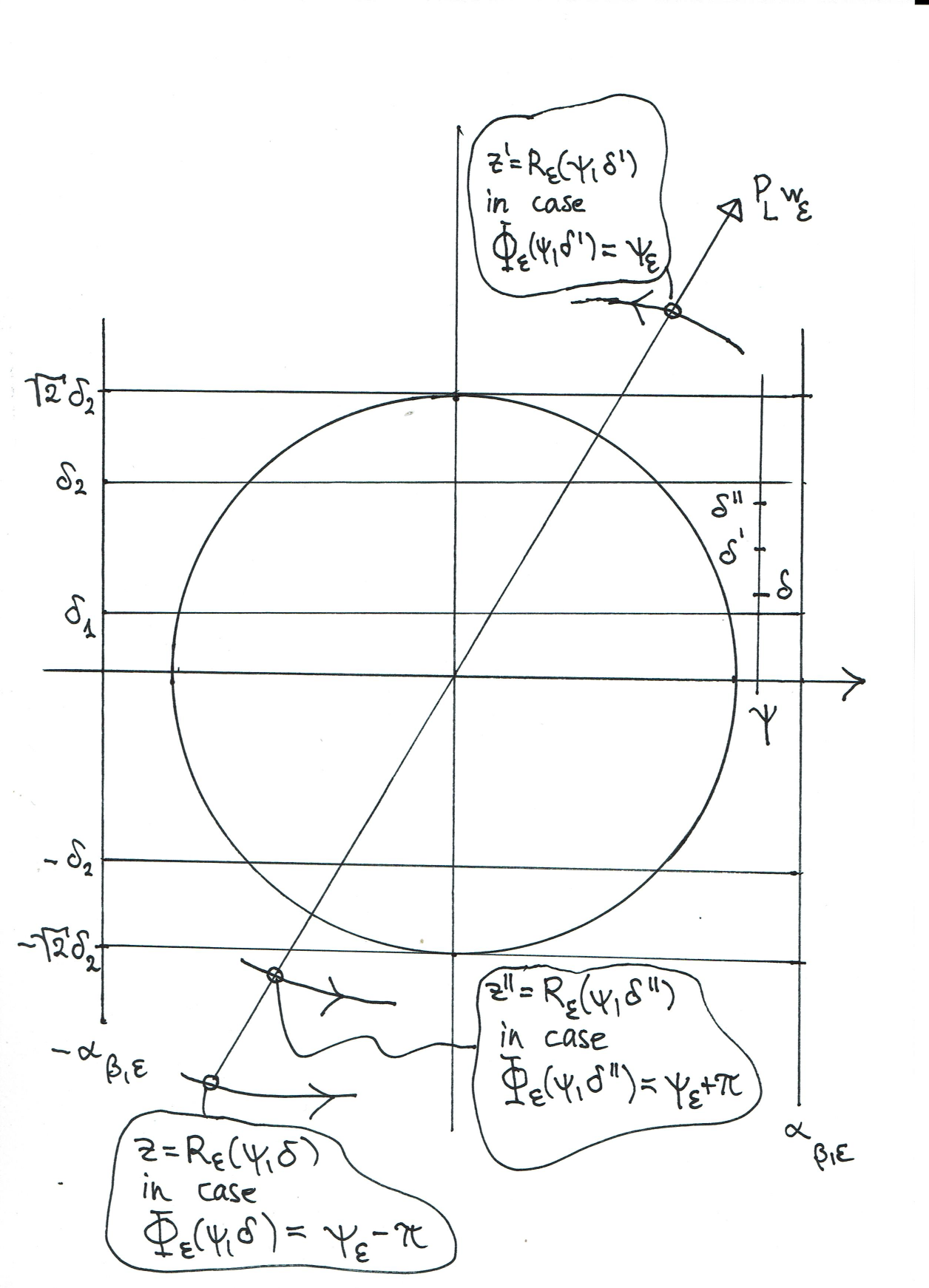}
	\caption{Positions of $R_{\epsilon}(\psi,\delta)$ depending on the angle $\Phi_{\epsilon}(\psi,\delta)$}
\end{figure}

{\bf Proof.} 1. On assertion (i). 

\medskip

1.1. Using (17) in combination with $\omega_{I,\epsilon}<\pi$ and $\delta_ {\beta,\epsilon}<\alpha_{\beta}<\delta_{I,\epsilon,1}<1$ we obtain from Proposition 6.4 that for $0<\delta_1<\delta_2\le\delta_{\beta,\epsilon}$ the rectangle 
$[-\alpha_{\beta},\alpha_{\beta}]\times[\delta_1,\delta_2]$ 
is contained in the domain of definition of the diffeomorphisms $\kappa_{\epsilon} P_{LE}(I_{\epsilon}(K_{\epsilon}(\cdot,\cdot)))$ and $R_{\epsilon}$.

\medskip

1.2. By Proposition  5.3 (ii), 
$$
x=\kappa_{\epsilon}P_{LE}I_{\epsilon}(K_{\epsilon}(\psi,\delta))
$$
satisfies
$$
|x|\ge\frac{1}{|\kappa_{\epsilon}^{-1}|}\delta^{\frac{-\sigma+\eta}{u-\eta}}\ge
\frac{1}{|\kappa_{\epsilon}^{-1}|}\delta_1^{\frac{-\sigma+\eta}{u-\eta}}
=\frac{1}{|\kappa_{\epsilon}^{-1}|}k_{\eta}^{\frac{-\rho+\eta}{u-\eta}}\delta_2^{c_{\eta}\frac{-\sigma+\eta}{u-\eta}}>2\sqrt{2}\delta_2.
$$
From Proposition 6.4 and Corollary 6.3 we know that $x$
is contained in the domain of the map $K_{\epsilon}^{-1}(E_{\epsilon}((\kappa_{\epsilon}\circ P_{LE})^{-1}(\cdot)))$. It follows that
$$
z=R_{\epsilon}(\psi,\delta)=K_{\epsilon}^{-1}(E_{\epsilon}(I_{\epsilon}(K_{\epsilon}(\psi,\delta))))=K_{\epsilon}^{-1}(E_{\epsilon}((\kappa_{\epsilon}\circ P_{LE})^{-1}(x))).
$$
From Corollary 6.3, $|z-x|\le\beta|x|$, and we obtain 
$$
|z|\ge|x|-\beta|x|=(1-\beta)|x|\ge\frac{1}{2}|x|>\sqrt{2}\delta_2.
$$

\medskip

2. On assertion (ii) for $(\psi,\delta)\in[-\alpha_{\beta,\epsilon},\alpha_{\beta,\epsilon}]\times[\delta_1,\delta_2]$ with $\Phi_{\epsilon}(\psi,\delta)=\psi_{\epsilon}$. Let $z=R_{\epsilon}(\psi,\delta)\in\mathbb{R}^2$.

\medskip

2.1. The image
$P_{LE}I_{\epsilon}(K_{\epsilon}(\psi,\delta))$ is a positive multiple of
$$
\left(\begin{array}{c}\cos(\Phi_{\epsilon}(\psi,\delta))\\ \sin(\Phi_{\epsilon}(\psi,\delta))\\0\end{array}\right)=\left(\begin{array}{c}\cos(\psi_{\epsilon})\\ \sin(\psi_{\epsilon})\\0\end{array}\right)\in(0,\infty)w_{\epsilon},
$$
hence 
$$
\kappa_{\epsilon} P_{LE}I_{\epsilon}(K_{\epsilon}(\psi,\delta))\in(0,\infty)\left(\begin{array}{c}0\\1\end{array}\right).
$$ 

2.2. We use the previous abbreviations $x$ and $z$ and have 
$|z|\ge(1-\beta)|x|>\sqrt{2}\delta_2$ from Part 1.
By Part 2.1, $x=\left(\begin{array}{c}0\\x_2\end{array}\right)$ with $x_2>0$.

\medskip

2.3. Proof of $|z|\le\sqrt{2}z_2$: We have $|x|=x_2$. From 
$x_2-z_2\le|x_2-z_2|\le|z-x|\le\beta|x|=\beta x_2$, $z_2\ge(1-\beta)x_2>0$. Also, from $x_1=0$, $|z_1|\le|x_1|+\beta|x|=\beta x_2$. It follows that
$$
|z|^2=z_1^2+z_2^2\le\beta^2x_2^2+z_2^2\le\frac{\beta^2}{(1-\beta)^2}z_2^2+z_2^2\le 2z_2^2.
$$

2.4. Consequently, $z_2=|z_2|\ge\frac{1}{\sqrt{2}}|z|>\delta_2$.

\medskip

3. The proofs of assertion (ii) in the two remaining cases are analogous, making use of the fact that in both cases we have that
$\kappa_{\epsilon} P_{LE}I_{\epsilon}(K_{\epsilon}(\psi,\delta))$ is a positive multiple of 
$\left(\begin{array}{c}0\\-1\end{array}\right)$. $\Box$

\medskip

The next result makes precise what was briefly announced at the begin of the section. The disjoint sets mentioned there will be given in terms of the angle $\Phi_{\epsilon}(\psi,\delta)$ corresponding to the value of the inner map $I_{\epsilon}(x)$, $x=K_{\epsilon}(\psi,\delta)$, only, not by the position of the value $R_{\epsilon}(\psi,\delta)$ of the full return map in coordinates in the left or right halfplane. Our choice of disjoint sets  circumvents a discussion how the latter, namely, positions of values $R_{\epsilon}(\psi,\delta)$ left or right of the vertical axis, are related to the more accessible angles $\Phi_{\epsilon}(\psi,\delta)$.

\begin{proposition}
Assume H$_3(\eta,\epsilon)$ and (16). Consider $\alpha_{\beta,\epsilon}$ and $\delta_{\beta,\epsilon}$ as in Proposition 7.3. Assume (18) for $\eta>0$, and (19) for $\delta_2$. Set $\delta_1=\delta_{1,\delta_2,\eta}$. Consider the disjoint sets
$$
M_0= \{(\psi,\delta)\in[-\alpha_{\beta,\epsilon},\alpha_{\beta,\epsilon}]\times[\delta_1,\delta_2]:\psi_{\epsilon}-\pi<\Phi_{\epsilon}(\psi,\delta)<\psi_{\epsilon}\}
$$
and
$$
M_1 = \{(\psi,\delta)\in[-\alpha_{\beta,\epsilon},\alpha_{\beta,\epsilon}]\times[\delta_1,\delta_2]:\psi_{\epsilon}<\Phi_{\epsilon}(\psi,\delta)<\psi_{\epsilon}+\pi\}.
$$
For every curve $c:[a,b]\to[-\alpha_{\beta,\epsilon},\alpha_{\beta,\epsilon}]\times[\delta_1,\delta_2]$ with $c(a)_2=\delta_1$ and $c(b)_2=\delta_2$ there exist $a_0<b_0<a_1<b_1$ in $[a,b]$ such that 
\begin{eqnarray*}
\mbox{on} &(a_0,b_0), & 
c(t)\in M_0\quad\mbox{and}\quad R_{\epsilon}(c(t))_2\in(\delta_1,\delta_2),\\
& \mbox{with} & R_{\epsilon}(c(a_0))_2=\delta_1\quad\mbox{and}\quad R_{\epsilon}(c(b_0))_2=\delta_2,\\
\mbox{while on} &(a_1,b_1), & c(t)\in M_1\quad\mbox{and}\quad R_{\epsilon}(c(t))_2\in(\delta_1,\delta_2),\\
& \mbox{with} & R_{\epsilon}(c(a_1))_2=\delta_2\quad\mbox{and}\quad R_{\epsilon}(c(b_1))_2=\delta_1.
\end{eqnarray*}
\end{proposition}

\begin{figure}
	\includegraphics[page=1,scale=0.7]{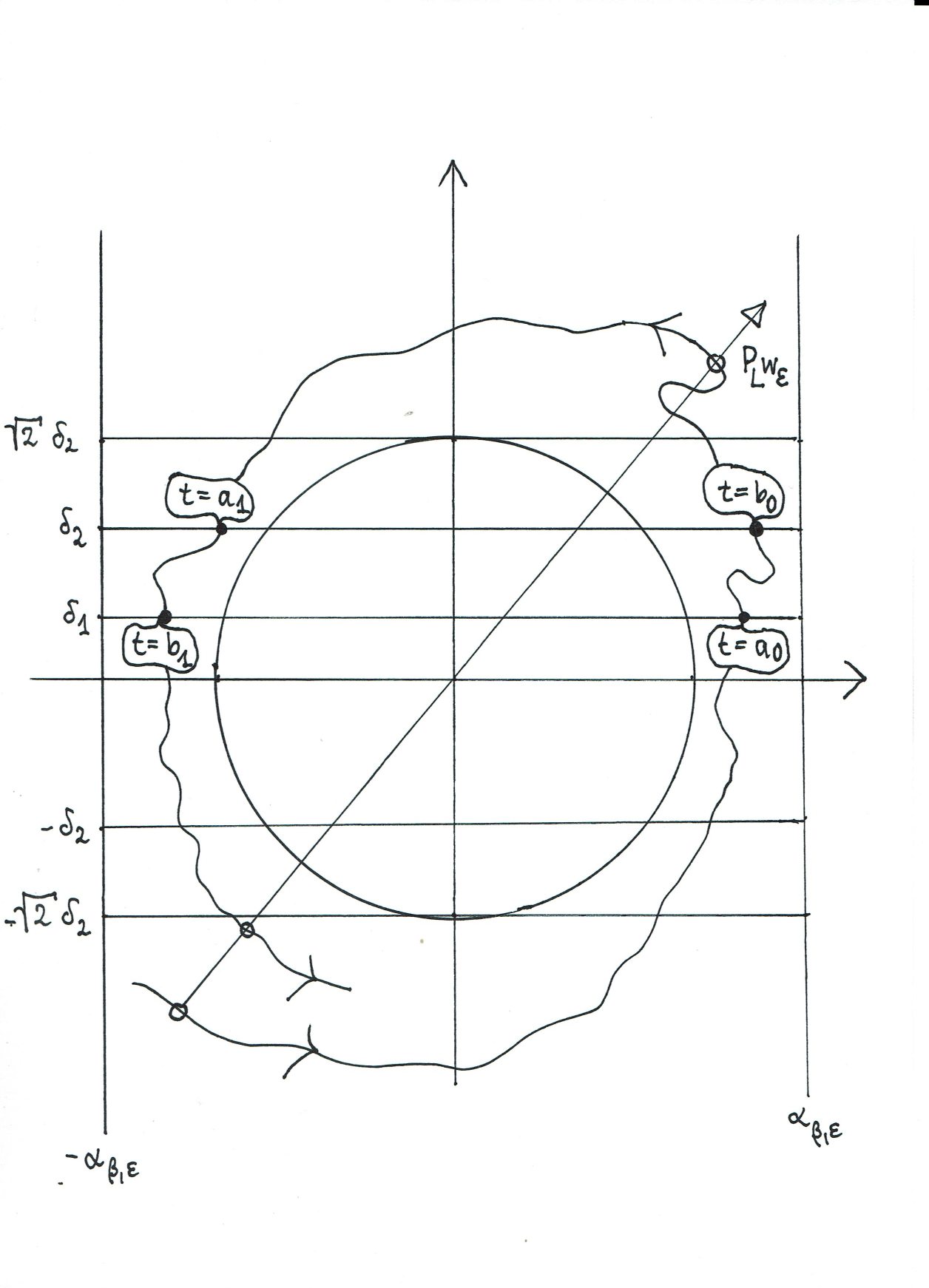}
	\caption{The values $R_{\epsilon}(c(t))$ for $a\le t\le b$}
\end{figure}

{\bf Proof.} 1. Proposition 7.2 yields $a_0'<b_0'\le a_1'<b_1'$ in $[a,b]$ such that
\begin{eqnarray*}
	\mbox{on} &(a_0',b_0'), & 
\Phi_{\epsilon}(c(t))\in(\psi_{\epsilon}-\pi,\psi_{\epsilon}),\\
& \mbox{with} & \Phi_{\epsilon}(c(a_0'))=\psi_{\epsilon}-\pi\quad\mbox{and}\quad \Phi_{\epsilon}(c(b_0'))=\psi_{\epsilon},\\	\mbox{and on} &(a_1',b_1'), & \Phi_{\epsilon}(c(t))\in(\psi_{\epsilon},\psi_{\epsilon}+\pi),\\
& \mbox{with} & \Phi_{\epsilon}(c(a_1'))=\psi_{\epsilon}\quad\mbox{and}\quad \Phi_{\epsilon}(c(b_1'))_2=\psi_{\epsilon}+\pi.
\end{eqnarray*}
From Proposition 7.3 (ii).
$$
R_{\epsilon}(c(a_0'))_2<-\delta_2,\quad R_{\epsilon}(c(b_0'))_2>\delta_2,\quad R_{\epsilon}(c(a_1'))_2>\delta_2,\quad R_{\epsilon}(c(b_1'))_2<-\delta_2.
$$
As in the proof of Proposition 7.2 one finds $a_0<b_0$ in $(a_0',b_0')$ and $a_1<b_1$ in $(a_1',b_1')$ with
\begin{eqnarray*}
R_{\epsilon}(c(a_0))_2=\delta_1 & \mbox{and} & R_{\epsilon}(c(b_0))_2=\delta_2,\quad\mbox{and}\quad R_{\epsilon}(c(t))\in(\delta_1,\delta_2)\quad\mbox{on}\quad(a_0,b_0),\\
R_{\epsilon}(c(a_1))_2=\delta_2 & \mbox{and} & R_{\epsilon}(c(b_1))_2=\delta_1,\quad\mbox{and}\quad R_{\epsilon}(c(t))\in(\delta_1,\delta_2)\quad\mbox{on}\quad(a_1,b_1).\\
\end{eqnarray*}
Observe that on $[a_0,b_0]\subset[a_0',b_0']$ we have $c(t)\in M_0$ while on $[a_1,b_1]\subset[a_1',b_1']$ we have $c(t)\in M_1$.
$\Box$

\section{Complicated dynamics}

For the results of this section we assume as in Proposition 7.4 that
H$_3(\eta,\epsilon)$ holds, and that $\beta$ satisfies (16). We also consider $\alpha_{\beta,\epsilon}$ and $\delta_{\beta,\epsilon}$ as in Proposition 7.3, and we assume (18) for $\eta>0$, and (19) for $\delta_2$. We set $\delta_1=\delta_{1,\delta_2,\eta}$ and consider the disjoint sets $M_0$ and $M_1$ from Proposition 7.4.

\begin{proposition}
For every sequence $(s_j)_0^{\infty}$ in $\{0,1\}$ there are forward trajectories $(x_j)_0^{\infty}$ of  $R_{\epsilon}$ with $x_j\in M_{s_j}$ and $\quad\delta_1\le R_{\epsilon}(x_j)_2\le\delta_2$ for all integers $\quad j\ge0$.
\end{proposition}

{\bf Proof.} 1. Let a sequence $(s_j)_0^{\infty}$ in $\{0,1\}$ be given. Choose a curve $c:[a,b]\to [-\alpha_{\beta,\epsilon} ,\alpha_{\beta,\epsilon}]\times[\delta_1,\delta_2]$ such that $c(t)_2\in(\delta_1,\delta_2)$ for $a<t<b$ and $c(a)_2=\delta_1$, $c(b)_2=\delta_2$, for example $c(t)=(0,t)$ for $a=\delta_1\le t\le\delta_2=b$.

\medskip

For integers $j\ge0$ we define recursively curves $c_j:[A_j,B_j]\to[-\alpha_{\beta,\epsilon},\alpha_{\beta,\epsilon}]\times[\delta_1,\delta_2]$ with decreasing domains in $[a,b]$ as follows.

\medskip

1.1. In order to define $c_0$ we apply Proposition 7.4 to the curve $c$ and obtain $a_0<b_0<a_1<b_1$ in $[a,b]$ with the properties stated in Proposition 7.4. In case $s_0=0$ we define $c_0$ by $A_0=a_0$, $B_0=b_0$, $c_0(t)=c(t)$ for $A_0\le t\le B_0$. Notice that $c_0(t)\in M_{s_0}$ for all $t\in(A_0,B_0)$,  $R_{\epsilon}(c_0(t))_2\in(\delta_1,\delta_2)$ on $(A_0,B_0)$,  $R_{\epsilon}(c_0(A_0))_2=\delta_1$, and $R_{\epsilon}(c_0(B_0))_2=\delta_2$. In case  $s_0=1$ we define $c_0$ by $A_0=a_1$, $B_0=b_1$, $c_0(t)=c(a_1+b_1-t)$ for $A_0\le t\le B_0$. Notice that also in this case $c_0(t)\in M_{s_0}$ for all $t\in(A_0,B_0)$,  $R_{\epsilon}(c_0(t))_2\in(\delta_1,\delta_2)$ on $(A_0,B_0)$,  and $R_{\epsilon}(c_0(A_0)_2)=R_{\epsilon}(c(a_1+b_1-a_1))_2=\delta_1$, $R_{\epsilon}(c_0(B_0))_2=R_{\epsilon}(c(a_1+b_1-b_1))_2=\delta_2$.

\medskip

1.2.  For an integer $j\ge0$ let a curve $c_j:[A_j,B_j]\to[-\alpha_{\beta,\epsilon},\alpha_{\beta,\epsilon}]\times[\delta_1,\delta_2]$ be given with  $c_j(t)\in M_{s_j}$ for all $t\in(A_j,B_j)$ and $R_{\epsilon}(c_j(t))_2\in(\delta_1,\delta_2)$ on $(A_j,B_j)$,  $R_{\epsilon}(c_j(A_j))_2=\delta_1$, $R_{\epsilon}(c_j(B_j))_2=\delta_2$. Proceeding as in Part 1.1, with the curve $[A_j,B_j]\ni t\mapsto R_{\epsilon}(c_j(t))\in[-\alpha_{\beta,\epsilon},\alpha_{\beta,\epsilon}]\times[\delta_1,\delta_2]$
in place of the former curve $c$, we obtain $A_{j+1}<B_{j+1}$ in $[A_j,B_j]$ and a curve $c_{j+1}:[A_{j+1},B_{j+1}]\to[-\alpha_{\beta,\epsilon},\alpha_{\beta,\epsilon}]\times[\delta_1,\delta_2]$ with  $c_{j+1}(t)\in M_{s_{j+1}}$ for all $t\in(A_{j+1},B_{j+1})$ and $R_{\epsilon}(c_{j+1}(t))_2\in(\delta_1,\delta_2)$ on $(A_{j+1},B_{j+1})$,  $R_{\epsilon}(c_{j+1}(A_{j+1}))_2=\delta_1$, $R_{\epsilon}(c_{j+1}(B_{j+1}))_2=\delta_2$.

\medskip

2. From $A_j\le A_{j+1}<B_{j+1}\le B_j$ for all integers $j\ge0$ we get $\cap_{j\ge0}[A_j,B_j]=[A,B]$ with $A=\lim_{j\to\infty}\le\lim_{j\to\infty}B_j=B$. By induction we obtain that for every $t\in[A,B]$ there is a forward trajectory $(x_j)_{j\ge0}$ of $R_{\epsilon}$
with $x_j=c_j(t)$, $A_j\le t\le B_j$, and $\delta_1\le R_{\epsilon}(x_j)_2\le\delta_2$ for all integers $j\ge0$. It remains to show that $x_j\in M_{s_j}$ for all integers $j\ge0$. Proof of this: Let an integer $j\ge0$ be given. From
$c_j(s)\in M_{s_j}$ for $A_j<s<B_j$ and $t\in[A_j,B_j]$ we get $x_j=c_j(t)\in cl M_{s_j}$. In case $s_j=0$ this yields $\psi_{\epsilon}-\pi\le\Phi_{\epsilon}(x_j)\le\psi_{\epsilon}$.
Assume $x_j\notin M_{s_j}$. Then $\Phi_{\epsilon}(x_j)\in\{\psi_{\epsilon}-\pi,\psi_{\epsilon}\}$,
and Proposition 7.3 (ii) yields $|R_{\epsilon}(x_j)|>\delta_2$, in contradiction to $\delta_1\le R_{\epsilon}(x_j)\le\delta_2$. In case $s_j=1$ the proof is analogous. $\Box$  

\medskip

The final result extends Proposition 8.1 to  entire trajectories.
  
\begin{theorem}
For every sequence $(s_j)_{j=-\infty}^{\infty}$ in $\{0,1\}$ there exist entire trajectories $(x_j)_{j=-\infty}^{\infty}$ of $R_{\epsilon}$ with
$x_j\in M_{s_j}$ for all integers $\,j$.
\end{theorem}

{\bf Proof.} 1. Let $(s_j)_{j=-\infty}^{\infty}$ in $\{0,1\}$ be given. Proposition 8.1 guarantees that for each integer $k$ there is a forward trajectory $(y_{k,j})_{j=0}^{\infty}$ of $R_{\epsilon}$ so that for all integers $j\ge0$,
$$
y_{k,j}\in M_{s_{j-k}}\quad\mbox{and}\quad\delta_1\le R_{\epsilon}(y_{k,j})_2\le\delta_2.
$$ 
For integers $k\ge-j$ we define
$$
z_{k,j}=y_{k,j+k},
$$
so that 
\begin{eqnarray*}
z_{k,j}=y_{k,j+k} & \in &  M_{s_{j+k-k}}=M_{s_j},\\ 
z_{k,j+1}=y_{k,j+1+k}=R_{\epsilon}(y_{k,j+k}) & = & R_{\epsilon}(z_{k,j}),\\
R_{\epsilon}(z_{k,j})_2=R_{\epsilon}(y_{k,j+k})_2 & \in & [\delta_1,\delta_2].
\end{eqnarray*}

1.1. Choice of subsequences for integers $J\ge0$.

\medskip

1.1.1. The case $J=0$: For all integers $k\ge0$, $z_{k,0}=y_{k,k}\in M_{s_0}$. As $[-\alpha_{\beta,\epsilon},\alpha_{\beta,\epsilon}]\times[\delta_1,\delta_2]$  is compact there exists a subsequence which converges to some $y_0\in\,cl\,M_{s_0}$. This subsequence is given by  
$z_{\kappa_0(k),0}$, $k\in\mathbb{N}_0$, with $\kappa_0:\mathbb{N}_0\to\mathbb{N}_0$ strictly
increasing.

\medskip

1.1.2. The case $J=1$: Consecutively choosing two further convergent subsequences we find a strictly increasing map
$\kappa_1:\mathbb{N}_0\to\mathbb{N}_0$ so that
for $k\to\infty$,
$$
z_{\kappa_0\circ\kappa_1(k),-1}\to y_{-1}\in\,cl\,M_{s_{-1}}\quad\mbox{and}\quad z_{\kappa_0\circ\kappa_1(k),1}\to y_1\in\,cl\,M_{s_1}.
$$

1.1.3. The general case $J\in\mathbb{N}_0$: Consecutively choosing further convergent subsequences analogously to Part 1.1.2 we obtain a strictly increasing map $\kappa_J:\mathbb{N}_0\to\mathbb{N}_0$
with $\kappa_J(0)\ge J$ so that for each $j\in\{-J,\ldots,J\}$ the sequence
$$
(z_{\kappa_0\circ\ldots\circ\kappa_J(k),j})_{k=0}^{\infty}
$$
converges for $k\to\infty$ to some $y_j\in\,cl\,M_{s_j}$.

\medskip

(Notice that for all integers $k\ge0$ we have 
$$
-J\ge-\kappa_J(0)\ge-\kappa_J(k)\ge-(\kappa_1\circ\ldots\circ\kappa_J)(k).)
$$

2. The {\it diagonal sequence} $K:\mathbb{N}_0\to\mathbb{N}_0$ defined by $K(J)=(\kappa_1\circ\ldots\circ\kappa_J)(J)$ is strictly increasing since for every $J\in\mathbb{N}$ we have
$$
K(J+1)=(\kappa_1\circ\ldots\circ\kappa_J)(\kappa_{J+1}(J+1))>
(\kappa_1\circ\ldots\circ\kappa_J)(\kappa_{J+1}(J))\ge
(\kappa_1\circ\ldots\circ\kappa_J)(J)=K(J)
$$
due to strict monotonicity of all maps involved.

\medskip

3. Let an integer $j$ be given and set $J=|j|$. In order to show that 
$$
(z_{K(k),j})_{k=J+1}^{\infty}\quad\mbox{is a subsequence of}\quad
(z_{(\kappa_1\circ\ldots\circ\kappa_J)(k),j})_{k=J+1}^{\infty}
$$
consider $\lambda:\{k\in\mathbb{N}_0:k>J\}\to\mathbb{N}_0$ given by
$$
\lambda(k)=(\kappa_{J+1}\circ\ldots\circ\kappa_k)(k).
$$
As in Part 2 one sees that the map $\lambda$ is strictly increasing, and
for every integer $k>J$,
$$
K(k)=(\kappa_1\circ\ldots\circ\kappa_k)(k)=((\kappa_1\circ\ldots\circ\kappa_J)\circ(\kappa_{J+1}\circ\ldots\circ\kappa_k))(k)=((\kappa_1\circ\ldots\circ\kappa_J)\circ\lambda)(k).
$$
It follows that $(z_{K(k),j})_{k=J+1}^{\infty}$ is a subsequence of $
(z_{(\kappa_1\circ\ldots\circ\kappa_{|j|})(k),j})_{k=J+1}^{\infty}$.

\medskip

4. We show that $(y_j)_{j=-\infty}^{\infty}$ is an entire trajectory of $R_{\epsilon}$. Let an integer $j$ be given and set $J=|j|$. From Part 3 in combination with Part 1.1.3 we get that
$(z_{K(k),j})_{k=J+1}^{\infty}$ converges to $y_j\in[-\alpha_{\beta,\epsilon},\alpha_{\beta,\epsilon}]\times[\delta_1,\delta_2]$ and that
$(z_{K(k),j+1})_{k=J+2}^{\infty}$ converges to $y_{j+1}$.
According to Part 1,
$z_{k,j+1}=R_{\epsilon}(z_{k,j})$ for all integers $k\ge-j$. For integers $k>J=|j|$ we have $j+1>j\ge-k\ge-K(k)$,
and the preceding statement yields
$$
z_{K(k),j+1}=R_{\epsilon}(z_{K(k),j}).
$$
It follows that
$$
y_{j+1}=\lim_{J+2\le k\to\infty}z_{K(k),j+1}=\lim_{J+1\le k\to\infty}R_{\epsilon}(z_{K(k),j})=R_{\epsilon}(y_j).
$$

\medskip

5. Proof of $y_j\in M_{s_j}$ for all integers $j$. Let an integer $j$ be given. We have $y_j\in\,cl\,M_{s_j}$. In case $s_j=0$ this yields $\psi_{\epsilon}-\pi\le\Phi_{\epsilon}(y_j)\le\psi_{\epsilon}$. Therefore the assumption $y_j\notin M_{s_j}$ results in  $\Phi_{\epsilon}(y_j)\in\{\psi_{\epsilon}-\pi,\psi_{\epsilon}\}$, which according to Proposition 7.3 (ii) means $|R_{\epsilon}(y_j)_2|>\delta_2$, in contradiction to $\delta_1\le R_{\epsilon}(y_j)_2\le\delta_2$. The proof in case $s_j=1$ is analogous. $\Box$
   
\section{Appendix: How to achieve (A1) and (A2)} 

Consider Shilnikov's scenario according to Section 1, with a twice continuously differentiable vectorfield $V$ and a homoclinic solution $h$ of Eq. (1). 

\medskip

The form of the linearization at $0$ in property (A1) results from replacing the  vectorfield $V$ by $V_{{\mathcal I}}:y\mapsto {\mathcal I}V({\mathcal I}^{-1}(y))$ with an isomorphism ${\mathcal I}:\mathbb{R}^3\to\mathbb{R}^3$ so that
$$
DV_{{\mathcal I}}(0)x={\mathcal I}DV(0){\mathcal I}^{-1}x=Ax
$$
for all $x\in\mathbb{R}^3$. Obviously, $V_{{\mathcal I}}(0)=0$. The  equation
$$
y'(t)=V_{{\mathcal I}}(y(t))
$$
is equivalent to Eq. (1) since $y$ is a solution if and only if $x={\mathcal I}^{-1}\circ y$ solves Eq. (1). Due to ${\mathcal I}(0)=0$ the  solution $h_{{\mathcal I}}={\mathcal I}\circ h$ satisfies
$\lim_{|t|\to\infty}h_{{\mathcal I}}(t)=0$. Obviously $h_{{\mathcal I}}(t)\neq0$ everywhere.
Let
$F_{{\mathcal I}}:\mathbb{R}^3\supset\,dom_{{\mathcal I}}\to\mathbb{R}^3$ denote the flow generated by the previous differential equation. $V_{{\mathcal I}}$ and $F_{{\mathcal I}}$ are twice continuously differentiable.

\medskip

Properties (3)-(6) can be achieved by a diffeomorphism ${\mathcal S}:\mathbb{R}^3\to\mathbb{R}^3$ which preserves what has been obtained by means of the isomorphism ${\mathcal I}$. A suitable diffeomorphism ${\mathcal S}$ can be found as follows.
Choose a neighbourhood $N$ of the origin on which $V_{{\mathcal I}}$ is close to its linearization  $DV_{{\mathcal I}}(0)$, say, $|V_{{\mathcal I}}(x)-Ax|<\epsilon|x|$ on $N$. Upon that change $V_{{\mathcal I}}$ outside an open neighourhood $N'\subset N$ of the origin to a twice continuously differentiable vectorfield $V':\mathbb{R}^3\to\mathbb{R}^3$
with $V'(x)=DV_{{\mathcal I}}(0)x=Ax$ outside 
$N$ and $|V'(x)-Ax|<\epsilon|x|$ on $N$. For $\epsilon>0$ sufficiently small the flow of $V'$ has global stable and unstable manifolds $W^s(0)\subset\mathbb{R}^3$ and $W^u(0)\subset\mathbb{R}^3$ of the stationary point $0$ which are invariant under that flow and have the form 
$$
W^s(0)=\{y+w^s(y):y\in L\},\quad W^u(0)=\{w^u(u)+u:u\in U\}
$$ 
with twice continuously differentiable maps $w^s:L\to U$ and $w^u:U\to L$ satisfying $w^s(0)=0,\,w^u(0)=0,\, Dw^s(0)=0,\, Dw^u(0)=0$. The intersections $W^s(0)\cap N'$ and $W^u(0)\cap N'$ are invariant in $N'$ under the flow $F_{{\mathcal I}}$ and have the property that solutions $y:(-\infty,t_0)\to N'$, $t_0\le\infty$, of the equation $y'(t)=V_{{\mathcal I}}(y(t))$
satisfy $y(t)\in W^u(0)$ on some unbounded interval $(-\infty,t_u)$, $t_u\le t_0$, while solutions $y:(t_0,\infty)\to N'$, $-\infty\le t_0$, of $y'(t)=V_{{\mathcal I}}(y(t))$
satisfy $y(t)\in W^s(0)$ on some unbounded interval $(t_s,\infty)$, $t_0\le t_s$. In particular, $h_{{\mathcal I}}(t)\in W^u(0)\cap N'$ on some interval $(-\infty,t_u)$ and $h_{{\mathcal I}}(t)\in W^s(0)\cap N'$ on some interval $(t_s,\infty)$.

\medskip

The map ${\mathcal S}:\mathbb{R}^3\to\mathbb{R}^3$ defined by ${\mathcal S}(y+u)=y-w^u(u)+u-w^s(y),\,y\in L,\,u\in U$ is a diffeomorphism with ${\mathcal S}(0)=0$ and $D{\mathcal S}(0)=id$ which transforms $W^s(0)$ onto $L$ and $W^u(0)$ onto $U$. $L$ and $U$ are invariant under the transformed flow 
$F_{{\mathcal S}}:\mathbb{R}^3\supset\,dom_{{\mathcal S}}\to\mathbb{R}^3$ which is defined by  
$(t,y)\in dom_{{\mathcal S}}$ if and only of $(t,{\mathcal S}^{-1}(y))\in dom_{{\mathcal I}}$, and in this case, $F_{{\mathcal S}}(t,y)={\mathcal S}(F_{{\mathcal I}}(t,{\mathcal S}^{-1}(y)))$. 
For the once continuously differentiable vectorfield
$$
V_{{\mathcal S}}:\mathbb{R}^3\ni z\mapsto D{\mathcal S}({\mathcal S}^{-1}(z))V_{{\mathcal I}}({\mathcal S}^{-1}(z))
\in\mathbb{R}^3
$$
we have 
$F_{{\mathcal S}}(t,y)=z(t)$ with the maximal solution
$z:I_y\to\mathbb{R}^3$ of the equation
$$
z'(t)=V_{{\mathcal S}}(z(t))
$$
with initial value $z(0)=y$. The preceding differential equation is equivalent to $y'(t)=V_{{\mathcal I}}(y(t))$ since $y$ is a solution of the latter if and only if  
$z={\mathcal S}\circ y$  is a solution of the former. Obviously, $V_{{\mathcal S}}(0)=0$, and (A1) holds for $V_{{\mathcal S}}$.

\medskip

The set ${\mathcal S}(N')$ is an open neighbourhood of the origin, and $L'={\mathcal S}(W^s(0))\cap {\mathcal S}(N')\subset L$ and $U'={\mathcal S}(W^u(0))\cap {\mathcal S}(N')\subset U$ are open neighbourhoods of the origin in $L$ and in $U$, respectively, and they are both invariant in ${\mathcal S}(N')$ under the transformed flow $F_{\mathcal S}$. By means of the relation 
$V_{{\mathcal S}}(z)=D_1F_{\mathcal S}(0,z)1$ for every $z\in\mathbb{R}^3$ one finds $V_{{\mathcal S}}(L')\subset L$ and $V_{{\mathcal S}}(U')\subset U$.

\medskip

Moreover, ${\mathcal S}(N'), L',$ and $U'$ have the property that solutions $z:(-\infty,t_0)\to {\mathcal S}(N')$, $t\le\infty$, of the equation $z'(t)=V_{{\mathcal S}}(z(t))$
satisfy $z(t)\in U'$ on some unbounded interval $(-\infty,t_u)$, $t_u\le t_0$, while solutions $z:(t_0,\infty)\to {\mathcal S}(N')$, $-\infty\le t_0$, of $z'(t)=V_{{\mathcal S}}(z(t))$ satisfy $z(t)\in L'$ on some unbounded interval $(t_s,\infty)$, $t_0\le t_s$. For $h_{{\mathcal S}}={\mathcal S}\circ h_{{\mathcal I}}$ we obtain from ${\mathcal S}(0)=0$ that $h_{{\mathcal S}}(t)\neq0$ everywhere and $lim_{|t|\to\infty}h_{{\mathcal S}}(t)=0$, and furthermore $h_{{\mathcal S}}(t)\in U'$ on some interval $(-\infty,t_u)$ and $h_{{\mathcal S}}(t)\in L'$ on some interval $(t_s,\infty)$.

\medskip

From the previous considerations it becomes obvious that
$V_{{\mathcal S}}$ has properties (A1) and (A2), with the homoclinic solution $h_{{\mathcal S}}$. Notice that due to (A1) also hypothesis (H) from Shilnikov's scenario is satisfied.

\medskip

Incidentally, notice that the flow $F_{\mathcal S}$ is as smooth as $V_{{\mathcal I}}$ and $F_{{\mathcal I}}$ (hence as smooth as $V$ and $F$) while  
$V_{{\mathcal S}}$ will in general not be better than once continuously differentiable. 

\medskip

Let us mention another possibility how to achieve property (A2) by a transformation. One could begin with local stable and unstable manifolds $W^s_{loc}(0)$ and $W^u_{loc}(0)$  of $F_{{\mathcal I}}$, which are given by maps $w^s_{loc}:L'\to U$ and $w^u_{loc}:U'\to L$ on neighbourhoods $L',\,U'$ of the origin in $L$ and in $U$, respectively. 
Restrictions of $w^s_{loc}$ and $w^u_{loc}$ to smaller neighbourhoods $L''\subset L'$ and $U''\subset U'$ would have twice continuously differentiable extensions to $L$ and $U$, respectively, which are zero outside $L'$ and $U'$.
Using the global maps $L\to U$ and $U\to L$ one would obtain a suitable diffeomorphism as above. A disadvantage of this second approach is that it does not generalize to semiflows in arbitrary Banach spaces, due to lack of a smooth norm which would be required for the construction of the extensions of $w^s_{loc}|L''$ and $w^u_{loc}|U''$  but
is not available in cases of infinite dimension.

\end{document}